\documentclass[twoside,12pt]{article}
\usepackage{amssymb,amsmath,bm,euscript}
\usepackage{graphicx,color,caption}
\usepackage{ifpdf}

\usepackage[colorlinks=true]{hyperref}

\usepackage{braket}

\usepackage{algorithm, algorithmic}

\newtheorem{proposition}{Proposition}[section]
\newtheorem{theorem}[proposition]{Theorem}
\newtheorem{lemma}[proposition]{Lemma}

\newtheorem{definition}[proposition]{Definition}

\newtheorem{example}[proposition]{Example}

\newcommand{\qed}{\hphantom{.}\hfill $\Box$\medbreak}

\def\F{\mathcal{F}}

\def\T{\mathcal{T}}

\def\A{{\mathcal{A}}}
\def\B{\mathcal{B}}
\def\C{\mathcal{C}}
\def\D{\mathcal{D}}

\def\X{{\mathcal{X}}}

\def\dd{{\bf d}}

\def\x{{\bf x}}
\def\y{{\bf y}}

\def\0{{\bf 0}}

\def\st{\mathrm{s.t.~}}

\title{\bf{Triple Decomposition and Tensor Recovery of Third Order Tensors}}
\author{ \hspace{1mm} Liqun Qi\thanks{
Department of Applied
    Mathematics, The Hong Kong Polytechnic University, Hung Hom,
    Kowloon, Hong Kong, China; ({\tt liqun.qi@polyu.edu.hk}).},
 \ \
 Yannan Chen\thanks{School of Mathematical Sciences, South China    Normal University, Guangzhou, China; ({\tt ynchen@scnu.edu.cn}).  This author was supported by the National Natural Science Foundation of China (11771405).},
 \ \
  Mayank Bakshi\thanks{Future Network Theory Lab, 2012 Labs
Huawei Tech. Investment Co., Ltd, Shatin, New Territory, Hong Kong, China; ({\tt mayank.bakshi@huawei.com}).}
 \ and \
 Xinzhen Zhang\thanks{School of Mathematics, Tianjin University, Tianjin 300354 China; ({\tt xzzhang@tju.edu.cn}). This author's work was supported by NSFC (Grant No.  11871369). }
}

\begin{document}
\date{\today}
\maketitle

\begin{abstract}
In this paper, we introduce a new tensor decomposition for third order tensors, which decomposes a third order tensor to three third order low rank tensors in a balanced way.   We call such a decomposition the triple decomposition, and the corresponding rank the triple rank.  For a third order tensor, its CP decomposition can be regarded as a special case of its triple decomposition.  The triple rank of a third order tensor is not greater than the middle value of the Tucker rank, and is strictly less than the middle value of the Tucker rank for an essential class of examples. These indicate that practical data can be approximated by low rank triple decomposition as long as it can be approximated by low rank CP or Tucker decomposition.    This theoretical discovery is confirmed numerically.  Numerical tests show that third order tensor data from practical applications such as internet traffic and video image are of low triple ranks.   A tensor recovery method based on low rank triple decomposition is proposed. Its convergence and convergence rate are established.     Numerical experiments confirm the efficiency of this method.


\end{abstract}


\section{Introduction}

Higher Order tensors have found many applications in recent years. Third order tensors are the most useful higher order tensors in applications~\cite{ADKM11, KBHH13, TYFWR13, XWWXWZ16, YHHH16, ZSKA18, ZA17, ZEAHK14, ZZXC15, ZLLZ18}. Tensor decomposition has emerged as a valuable tool for analyzing and computing with such tensors \cite{KB09}. For example, a key idea behind tensor recovery algorithms is that many practical datasets are highly structured in the sense that the corresponding tensors can be approximately represented through a low rank decomposition.


Two most well-known tensor decompositions are the CANDECOMP/PARAFAC (CP) decomposition and the Tucker decomposition \cite{KB09}.  Their corresponding ranks are called CP rank and Tucker rank \cite{JYZ17} respectively.  In the next section, we will review their definitions.

Suppose that we have a third order tensor $\X \in \Re^{n_1 \times n_2 \times n_3}$, where $n_1, n_2$ and $n_3$ are positive integers.   The CP rank of $\X$ may be higher than max$\{ n_1, n_2, n_3 \}$.  For example, the CP rank of a $9 \times 9 \times 9$ tensor given by Kruskal is between $18$ and $23$.   See \cite{KB09}.  It is known \cite{KB09} that an upper bound of the CP rank is min$\{ n_1n_2, n_1n_3, n_2n_3 \}$.

The Tucker decomposition decomposes $\X$ into a core tensor $\D \in \Re^{r_1\times r_2 \times r_3}$ multiplied by three factor matrices $U \in \Re^{n_1 \times r_1}$, $V \in \Re^{n_2 \times r_2}$ and $W \in \Re^{n_3 \times r_3}$ along three modes, i.e.,
$$\X = \D \times_1 U \times_2 V \times_3 W.$$
The minimum possible values of $r_1, r_2$ and $r_3$ are called the Tucker rank of $\X$ \cite{KB09}.   Then $r_i \le n_i$ for $i = 1, 2, 3$.   Thus, the Tucker rank is relatively smaller.

In this paper, we introduce a new tensor decomposition for third order tensors, which decomposes a third order tensor to a product of three third order low rank tensors in a balanced way.   We call such a decomposition the triple decomposition, and the corresponding rank the triple rank. For a third order tensor, its CP decomposition can be regarded as a special case of its triple decomposition.  The triple rank of a third order tensor is not greater than the middle value of the Tucker rank, and is strictly less than the middle value of the Tucker rank for an essential class of examples. These indicate that practical data can be approximated by low rank triple decomposition as long as it can be approximated by low rank CP or Tucker decomposition.   This theoretical discovery is  confirmed numerically.  Numerical tests show that third order tensor data from practical applications such as internet traffic and video image are of low triple ranks.   A tensor recovery method based on low rank triple decomposition is proposed. Its convergence and convergence rate are established.     Numerical experiments confirm the efficiency of this method.

The rest of this paper is distributed as follows.     Preliminary knowledge on CP decomposition, Tucker decomposition, and related tensor ranks is presented in the next section.  In Section 3, we introduce triple decomposition and triple rank, prove the above key properties and some other theoretical properties.  In particular, we show that the triple rank of a third order tensor is not greater than the triple rank of the Tucker core of that third order tensor.  If the factor matrices of the Tucker decomposition of that third order tensor are of full column rank, than the triple ranks of that third order tensor and its Tucker core are equal.  We present an algorithm to check if a given third order tensor can be approximated by a third order tensor of low triple rank such that the relative error is reasonably small, and make convergence analysis for this algorithm in Section 4.    In Section 5, we show that practical data of third order tensors from  internet traffic and video image are of low triple ranks.   A tensor recovery method is proposed in Section 6, based on such low rank triple decomposition.  Its convergence and convergence rate are also established in that section.   Numerical comparisons of our method with tensor recovery based upon CP and Tucker decompositions are presented in Section 7.   Some concluding remarks are made in Section 8.


\section{CP Decomposition, Tucker Decomposition and Related Tensor Ranks}

We use small letters to denote scalars, small bold letters to denote vectors, capital letters to denote matrices, and calligraphic letters to denote tensors.   In this paper, we only study third order tensors.

Perhaps, the most well-known tensor decomposition is CP decomposition \cite{KB09}.    Its corresponding tensor rank is called the CP rank.

\begin{definition} \label{d1}
Suppose that $\X = (x_{ijt}) \in \Re^{n_1 \times n_2 \times n_3}$.   Let $A = (a_{ip}) \in \Re^{n_1 \times r}$, $B = (b_{jp}) \in \Re^{n_2 \times r}$ and $C = (c_{tp}) \in \Re^{n_3 \times r}$.  Here $n_1, n_2, n_3, r$ are positive integers.   If
\begin{equation} \label{eq2.1}
x_{ijt} = \sum_{p=1}^r a_{ip}b_{jp}c_{tp}
\end{equation}
for $i=1, \cdots, n_1$, $j = 1, \cdots, n_2$ and $t = 1, \cdots, n_3$, then $\X$ has a CP decomposition $\X = [[A, B, C]]$.   The smallest integer $r$ such that (\ref{eq2.1}) holds is called the CP rank of $\X$, and denoted as CPRank$(\X) = r$.
\end{definition}

As shown in \cite{KB09}, CPRank$(\X) \le \min \{n_1n_2, n_1n_3, n_2n_3 \}$.   A tensor recovery method via CP decomposition can be found in \cite{ADKM11}.

Another well-known tensor decomposition is Tucker decomposition \cite{KB09}.    Its corresponding tensor rank is called the Tucker rank.     Higher order SVD (HOSVD) decomposition \cite{DDV00} can be regarded as a special variant of Tucker decomposition.

\begin{definition} \label{d2}
Suppose that $\X = (x_{ijt}) \in \Re^{n_1 \times n_2 \times n_3}$, where $n_1, n_2$ and $n_3$ are positive integers.   We may unfold $\X$ to a matrix $X_{(1)} = (x_{i, jt}) \in \Re^{n_1 \times n_2n_3}$, or a matrix $X_{(2)} = (x_{j, it}) \in \Re^{n_2 \times n_1n_3}$, or a matrix $X_{(3)} = (x_{t, ij}) \in \Re^{n_3 \times n_1n_2}$.   Denote the matrix ranks of $X_{(1)}, X_{(2)}$ and $X_{(3)}$ as $r_1, r_2$ and $r_3$, respectively.  Then the triplet $(r_1, r_2, r_3)$ is called the Tucker rank of $\X$, and is denoted as TucRank$(\X) = (r_1, r_2, r_3)$ with TucRank$(\X)_i = r_i$ for $i=1, 2, 3$.
\end{definition}

The CP rank and Tucker rank are called the rank and $n$-rank in some papers \cite{KB09}.  Here, we follow \cite{JYZ17} to distinguish them from other tensor ranks.

\begin{definition} \label{d3}
Suppose that $\X = (x_{ijt}) \in \Re^{n_1 \times n_2 \times n_3}$.   Let $U = (u_{ip}) \in \Re^{n_1 \times r_1}$, $V = (v_{jq}) \in \Re^{n_2 \times r_2}$, $W = (w_{ts}) \in \Re^{n_2 \times r_3}$ and $\D = (d_{pqs}) \in \Re^{r_1 \times r_2 \times r_3}$.  Here $n_1, n_2, n_3, r_1, r_2, r_3$ are positive integers.   If
\begin{equation} \label{eq2.2}
x_{ijt} = \sum_{p=1}^{r_1} \sum_{q=1}^{r_2} \sum_{s=1}^{r_3} u_{ip}v_{jq}c_{ts}w_{pqs}
\end{equation}
for $i=1, \cdots, n_1$, $j = 1, \cdots, n_2$ and $t = 1, \cdots, n_3$, then $\X$ has a Tucker decomposition $\X = [[\D; U, V, W]]$.   The matrices $U, V, W$ are called factor matrices of the Tucker decomposition, and the tensor $\D$ is called the Tucker core.   We may also denote the Tucker decomposition as
\begin{equation} \label{eq2.3}
\X = \D \times_1 U \times_2 V \times_3 W.
\end{equation}
\end{definition}

The Tucker ranks $r_1, r_2, r_3$ of $\X$ are the smallest integers such that (\ref{eq2.2}) holds \cite{KB09}.   Nonnegative tensor recovery methods via Tucker decomposition can be found in \cite{XY13, CSZZC19}.

\section{Triple Decomposition, Triple Rank and Their Properties}

Let $\X = (x_{ijt}) \in \Re^{n_1 \times n_2 \times n_3}$.  As in \cite{KB09}, we use $\X(i, :, :)$ to denote the $i$-th horizontal slice, $\X(:, j, :)$ to denote the $j$-th lateral slice; $\X(:, :, t)$ to denote the $t$-th frontal slice.  We say that $\X$ is a third order horizontally square tensor if all of its horizontal slices are square, i.e., $n_2 = n_3$. Similarly, $\X$ is a third order laterally square tensor (resp. frontally square tensor) if all of its lateral slices (resp. frontal slices) are square, i.e., $n_1 = n_3$ (resp.  $n_1 = n_2$).

\begin{definition} Let $\X = (x_{ijt}) \in \Re^{n_1 \times n_2 \times n_3}$ be a non-zero tensor.   We say that $\X$ is the triple product of a third order horizontally square tensor $\A = (a_{iqs}) \in \Re^{n_1 \times r \times r}$, a third order laterally square tensor $\B = (b_{pjs}) \in \Re^{r \times n_2 \times r}$ and a third order frontally square tensor $\C = (c_{pqt}) \in \Re^{r \times r \times n_3}$, and denote
\begin{equation}  \label{eq9}
\X = \A\B\C,
\end{equation}
if for $i = 1, \cdots, n_1$, $j = 1, \cdots, n_2$ and $t = 1, \cdots, n_3$, we have
\begin{equation}  \label{e2}
x_{ijt} = \sum_{p, q, s=1}^r a_{iqs}b_{pjs}c_{pqt}.
\end{equation}
If
\begin{equation}  \label{eq11}
r \le {\rm mid} \{ n_1, n_2, n_3 \},
\end{equation}
then we call (\ref{eq9}) a low rank triple decomposition of $\X$.  See Figure 1 for a visualization.
\input{triplerank.TpX}

The smallest value of $r$ such that (\ref{e2}) holds is called the triple rank of $\X$, and is denoted as TriRank$(\X) = r$.  For a zero tensor, we define its triple rank as zero.
\end{definition}

Note that TriRank$(\X)$ is zero if and only if it is a zero tensor.   This is analogous to the matrix case.

\begin{theorem}
Low rank triple decomposition and triple ranks are well-defined.  A third order nonzero tensor $\X$ always has a low rank triple decomposition (\ref{eq9}), satisfying (\ref{eq11}).
\end{theorem}
{\bf Proof}    Without loss of generality, we may assume that we have a third order nonzero tensor $\X \in \Re^{n_1 \times n_2 \times n_3}$ and $n_1 \ge n_2 \ge n_3 \ge 1$.  Thus, mid$\{ n_1, n_2, n_3 \} = n_2$.   Let $r = n_2$. Let $\A\in\Re^{n_1\times r\times r}$, $\B\in\Re^{r\times n_2\times r}$, and $\C\in\Re^{r \times r\times n_3}$ be such that $a_{iqs} = x_{isq}$ if $q \le n_3$, $a_{iqs} = 0$ if $q > n_3$, $b_{pjs} = {\delta_{js} \over r}$, and $c_{pqt} = {\delta_{qt} \over r}$ for $i = 1, \cdots, n_1$, $j, p, q, s = 1, \cdots, n_2$, and $t = 1, \cdots, n_3$, where $\delta_{js}$ and $\delta_{qt}$ are the Kronecker symbol such that $\delta_{jj} = 1$ and $\delta_{js} = 0$ if $j \not = s$.    Then (\ref{e2}) holds for the above choices of $\A$, $\B$, and $\C$. Thus, the triple decomposition always exists with $r\leq n_2$. \qed


Note that one cannot change (\ref{eq11}) to
\begin{equation} \label{eq12}
r \le \min \{ n_1, n_2, n_3 \}.
\end{equation}
The above assertion can be seen through the following argument.  Let $n_1 = n_2 = 3$ and $n_3 = 1$.  Suppose that $\X$ is chosen to have $9$ independent entries.  If
(\ref{eq12}) is required, then with $r = 1$, the decomposition consists of $\A$, $\B$ and $\C$ that can have only a maximum of $7$ independent entries in total.  Thus, we cannot find $\A$, $\B$, and $\C$, satisfying (\ref{eq9}), (\ref{e2}) and (\ref{eq12}).


Suppose that $\X=\A\B\C$, where $\A=\mathcal{F}\times_1 \widetilde{A}$, $\mathcal{F}\in\Re^{r_1\times r\times r}$, $\widetilde{A}\in\Re^{n_1\times r_1}$,
$\B=\mathcal{G}\times_2 \widetilde{B}$, $\mathcal{G}\in\Re^{r\times r_2\times r}$, $\widetilde{B}\in\Re^{n_2\times r_2}$,
$\C=\mathcal{H}\times_3 \widetilde{C}$, $\mathcal{H}\in\Re^{r\times r\times r_3}$ and $\widetilde{C}\in\Re^{n_3\times r_3}$.
Then, we have
\begin{equation} \label{tripletucker}
x_{ijk} = \sum_{p, q, s=1}^r a_{iqs}b_{pjs}c_{pqk}
  = \sum_u^{r_1}\sum_v^{r_2}\sum_w^{r_3}\widetilde{A}_{iu}\widetilde{B}_{jv}\widetilde{C}_{kw}\underbrace{\sum_{p, q, s=1}^rF_{uqs}G_{pvs}H_{pqw}}_{\text{a core tensor }\mathcal{FGH}}.
\end{equation}
Thus, this is a formulation of the Tucker decomposition. In addition, if the core tensor $\mathcal{FGH}$ is a diagonal tensor, we get the CP decomposition.

We now study the relation between triple decomposition and CP decomposition.   We have the following theorem.

\begin{theorem} \label{p1}
Suppose that $\X = (x_{ijt}) \in \Re^{n_1 \times n_2 \times n_3}$.   Then we may regard its CP decomposition as a special case of its triple decomposition.  In particular, we have
$${\rm TriRank}(\X) \le {\rm CPRank}(\X) \le {\rm TriRank}(\X)^3.$$
\end{theorem}
{\bf Proof} Suppose that $\X=[[A,B,C]]$ with $A=(a_{ip})\in \Re^{n_1\times r}$, $B=(b_{jp})\in \Re^{n_2\times r}$ and $C\in \Re^{n_3\times r}$
is a CP decomposition. Denote $\A=(\bar a_{ipq})\in \Re^{n_1\times r\times r}$, $\B=(\bar b_{sjq})\in \Re^{r\times n_2\times r}$ and $\C=(c_{spt})\in
\Re^{r\times r\times n_3}$ with
\[\begin{array}{ll}&\bar a_{ipq}=\left\{\begin{array}{lll}&a_{ip}\quad&{\mbox if} ~~p=q,\\
&0, &{\mbox otherwise}.\end{array}\right. \quad
\bar b_{sjq}=\left\{\begin{array}{lll}& b_{jq}\quad&{\mbox if}~~ s=q,\\
&0, &{\mbox otherwise}.\end{array}\right.\\
&\bar c_{spt}=\left\{\begin{array}{lll}& c_{tp} \quad&{\mbox if}~~ s=p,\\
&0, &{\mbox otherwise}.\end{array}\right.
\end{array}\]
Then for all
$i=1,\dots, n_1$, $j=1,\dots, n_2$ and $t=1,\dots, n_3$, there holds
$$(\A\B\C)_{ijt}=\sum\limits_{s,p,q=1}^r \bar a_{ipq}\bar b_{sjq}\bar c_{spt}=\sum\limits_{p=1}^r a_{ip}b_{jp}c_{tp}=\X_{ijt}.$$ This means that $\X=\A\B\C$, i.e., we may regard its CP decomposition as a special case of its triple decomposition.  Furthermore, we have
TriRank$(\X)\le$ CPRank$(\X)$ from the definition of the triple rank.

On the other hand, suppose that $\X$ is of the form $x_{ijt}=\sum_{p,q,s=1}^{\bar r}a_{iqs}b_{pjs}c_{pqt}$. Then, $\X$ can be represented as a sum of $\bar r^3$ rank-one tensors.   Hence, the last inequality in the theorem holds  by setting $\bar r=TriRank(\X)$.
\qed



This theorem indicates that the triple rank is not greater than the CP rank.  As the CP rank may be greater than max$\{ n_1, n_2, n_3 \}$, while the triple rank is not greater than mid$\{ n_1, n_2, n_3 \}$, there is a good chance that
 the triple rank is strictly smaller than the CP rank.   By \cite{KB09}, Monte Carlo experiments reveal that the set of $2 \times 2 \times 2$ tensors are of CP rank three with probability $0.21$.  Since the triple rank is not greater than two in this case, with a substantial probability the triple rank is strictly less than the CP rank.

 Next, we study the relation between triple decomposition and Tucker decomposition.  We have the following theorem.

\begin{theorem} \label{kt}
Suppose that $\X = (x_{ijt}) \in \Re^{n_1 \times n_2 \times n_3}$ and $\X=\A\B\C$  with TriRank$(\X)=R$, $\A\in \Re^{n_1\times R\times R}, \B \in \Re^{R\times n_2\times R}$,  $\C\in \Re^{R\times R\times n_3}$. Furthermore, $$\X=\D\times_1 U\times_2 V\times_3W$$ is a Tucker decomposition of $X$ with $\D\in \Re^{r_1\times r_2\times r_3}$ and factor matrices $U\in \Re^{n_1\times r_1}, V\in \Re^{n_2\times r_2}, W\in \Re^{n_3\times r_3}$.
Then
\begin{equation} \label{eqq3.10}
{\rm TriRank}(\X) \le {\rm TriRank}(\D) \le {\rm mid}\{ r_1, r_2, r_3 \}.
\end{equation}

Furthermore, if TucRank$(\X)=(r_1,r_2,r_3)$, then we have
\begin{equation} \label{eqq3.11}
{\rm TriRank}(\X) = {\rm TriRank}(\D).
\end{equation}
Thus, we always have
\begin{equation} \label{eqq3.12}
{\rm TriRank}(\X) \le {\rm mid} \{ {\rm TucRank}(\X)_1, {\rm TucRank}(\X)_2, {\rm TucRank}(\X)_3 \}.
\end{equation}
\end{theorem}
{\bf Proof} For convenience of notation, let TriRank$(\D)=r$.   By (\ref{eq11}), we have the second inequality of (\ref{eqq3.10}).

We first show that $r\geq R$.
Assume that $\D=\bar\A\bar \B\bar \C$ with $\bar\A\in \Re^{r_1\times r\times r}$, $\bar \B\in \Re^{r\times r_2\times r}$ and $\bar \C\in \Re^{r\times r\times r_3}$.  Then
$$\X=(\bar\A \bar \B\bar\C)\times_1 U\times_2V\times_3W
=(\bar\A\times_1 U)(\bar \B\times_2V)(\bar\C\times_3W).
$$
Clearly, $\bar\A\times_1 U\in \Re^{n_1\times r\times r}$,  $\bar \B\times_2V\in \Re^{r\times n_2\times r}$ and
$\bar\C\times_3W\in \Re^{r\times r\times  n_3}$. Hence, $r\geq R$ from the definition of TriRank.   This proves the first inequality of (\ref{eqq3.10}).

Now we assume that TucRank$(\X)=(r_1,r_2,r_3)$, and show that $r\leq R$.  By TucRank$(\X)=(r_1,r_2,r_3)$, we know that factor matrices $U, V$ and $W$ are of full column rank.   Then $U^TU, V^TV$ and $W^TW$ are invertible.
From $\X=\D\times_1 U\times_2 V\times_3 W$, we have that
$$\begin{array}{rl}& \X\times_1 (U^TU)^{-1}U^T \times_2 (V^TV)^{-1}V^T\times_3 (W^TW)^{-1}W^T\\
=&(\D \times_1 U\times_2 V\times_3 W)\times_1 (U^TU)^{-1}U^T \times_2 (V^TV)^{-1}V^T \times_3(W^TW)^{-1}W^T\\
=&\D \times_1 (U^TU)^{-1}(U^TU)\times_2 (V^T V)^{-1}(V^TV)\times_3 (W^TW)^{-1}(W^TW)\\
=&\D\times_1 I_{r_1}\times_2 I_{r_2}\times_3 I_{r_3}=\D.
 \end{array}$$
Hence, it holds that
$$\begin{array}{rl}\D=&(\A\B\C)\times_1 (U^TU)^{-1}U^T\times_2 (V^TV)^{-1}V^T\times_3(W^TW)^{-1}W^T\\
=& (\A\times_1 (U^TU)^{-1}U^T)(\B\times_2(V^TV)^{-1}V^T)(\C\times_3(W^TW)^{-1}W^T). \end{array}$$
It is easy to see that $\A\times_1 (U^TU)^{-1}U^T\in \Re^{r_1\times R\times R}$, $\B\times_2(V^TV)^{-1}V^T\in \Re^{R\times r_2\times R}$ and
$\C\times (W^TW)^{-1}W^T \in \Re^{R\times R\times r_3}$. From definition of TriRank,  we have $r\leq R$.
Therefore $r=R$ and (\ref{eqq3.11}) holds.

Note that the condition that TucRank$(\X)=(r_1,r_2,r_3)$ always can be realized.   For example, in HOSVD \cite{DDV00}, all factor matrices are orthogonal, and we always have TucRank$(\X)=(r_1,r_2,r_3)$.   This shows that (\ref{eqq3.12}) always holds.
\qed

The condition that TucRank$(\X)=(r_1,r_2,r_3)$ holds if all factor matrices are of full column rank.  In \cite{JYZ17}, if the factor matrices of a Tucker decomposition are of full column rank, then that Tucker decomposition is called independent.

We now give an example that TriRank$(\X) <$ min$\{$ TucRank$(\X)_1$, TucRank$(\X)_2$, TucRank$(\X)_3 \}$.

\begin{example}
Let $n_1 = n_2 = n_3 = 4$ and $r = 2$. Consider $\A = (a_{iqs}) \in \Re^{4 \times 2 \times 2}$, $\B = (b_{pjs}) \in \Re^{2 \times 4 \times 2}$, and
$\C = (c_{pqt}) \in \Re^{2 \times 2 \times 4}$ such that $a_{111} = a_{212} = a_{321} = a_{422} = 1$ and $a_{iqs} = 0$ otherwise, $b_{111} = b_{122} = b_{231} = b_{242} = 1$ and $b_{pjs} = 0$ otherwise, and $c_{111} = c_{122} = c_{213} = c_{224} = 1$ and $c_{pqt} = 0$ otherwise.  Then TucRank$(\A)_1 =$ TucRank$(\B)_2 =$ TucRank$(\C)_3 = 4$.  Let $\X = \A\B\C$.   Then TriRank$(\X) \le 2$ and $\X \in \Re^{4 \times 4 \times 4}$. We have $x_{111} = x_{133} = x_{221} = x_{243} = x_{312} = x_{334} = x_{422} = x_{444} = 1$ and $x_{ijt} = 0$ otherwise.    We may easily check that
TucRank$(\X)_1 =$ TucRank$(\X)_2 =$ TucRank$(\X)_3 = 4$.   Thus, TriRank$(\X) \le 2  <$ TucRank$(\X)_1 =$ TucRank$(\X)_2 =$ TucRank$(\X)_3 = 4$.
\end{example}

Taking the conclusion of the above example further, the following probabilistic argument shows that, in fact, the triple rank is smaller than the smallest Tucker rank for an essential class of examples.  Let $n_1 = n_2 = n_3 = 4$ and $r = 2$, $\A = (a_{iqs}) \in \Re^{4 \times 2 \times 2}$, $\B = (b_{pjs}) \in \Re^{2 \times 4 \times 2}$ and $\C = (c_{pqt}) \in \Re^{2 \times 2 \times 4}$.   Then $A_{(1)}, B_{(2)}$ and $C_{(3)}$ are $4 \times 4$ matrices.   With probability one, these three matrices are nonsingular, i.e., TucRank$(\A)_1 =$
TucRank$(\B)_2 =$ TucRank$(\C)_3 = 4$.    Let $\X = \A\B\C$.   Then $\X \in
\Re^{4 \times 4 \times 4}$ and $X_{(1)}, X_{(2)}$ and $X_{(3)}$ are $4 \times 4$ matrices.  With probability one, these three matrices are also nonsingular, i.e., TucRank$(\X)_1 =$ TucRank$(\X)_2 =$ TucRank$(\X)_3 = 4$.   Then, with probability one, we have TriRank$(\X) \le 2  <$ TucRank$(\X)_1 =$ TucRank$(\X)_2 =$ TucRank$(\X)_3 = 4$.    This shows that there is a substantial chance that TriRank$(\X) <$ mid$\{$ TucRank$(\X)_1$, TucRank$(\X)_2$, TucRank$(\X)_3 \}$.

The above two theorems indicate that practical data can be approximated by low rank triple decomposition as long as it can be approximated by low rank CP or Tucker decomposition.
This theoretical discovery will be confirmed numerically in the later sections.

Next, we have the following proposition relating the triple rank to the Tucker rank.

\begin{proposition} \label{p2}
Suppose that $\X = (x_{ijk}) \in \Re^{n_1 \times n_2 \times n_3}$ and $\X=\A\B\C$ with $\A\in \Re^{n_1\times r\times r}, \B \in \Re^{r\times n_2\times r}$ and $\C\in \Re^{r\times r\times n_3}$. Then
\begin{equation} \label{eq3.16}
{\rm TucRank}(\X)_1 \le {\rm TucRank}(\A)_1 \le \left({\rm TriRank}(\A)\right)^2 \le \left({\rm TriRank}(\X)\right)^2,
\end{equation}
\begin{equation} \label{eq3.17}
{\rm TucRank}(\X)_2 \le {\rm TucRank}(\B)_2 \le \left({\rm TriRank}(\B)\right)^2 \le \left({\rm TriRank}(\X)\right)^2,
\end{equation}
and
\begin{equation} \label{eq3.18}
{\rm TucRank}(\X)_3 \le {\rm TucRank}(\C)_3 \le \left({\rm TriRank}(\C)\right)^2 \le \left({\rm TriRank}(\X)\right)^2.
\end{equation}
\end{proposition}
{\bf Proof}  Let  TriRank$(\X)=r$, TucRank$(\A)_1=r_1$, TucRank$(\B)_2=r_2$ and TuckRank$(\C)_3=r_3$. Let
$\A=\F \times_1 U\times_2U_2\times_3 U_3$ be a Tucker decomposition of $\A$  with core tensor $\F\in\Re^{r_1\times s_2\times s_3} $ and factor matrices $U\in \Re^{n_1\times r_1}, U_2\in \Re^{r\times s_2}, U_3\in \Re^{r\times s_3}$. Denote $\bar \A=\F\times_2U_2\times_3 U_3\in \Re^{r_1\times r\times r}$. Then $\A=(\F\times_2U_2\times_3 U_3)\times_1 U=\bar\A \times_1 U$.

Similarly,
there exist $\bar\B\in \Re^{r\times r_2\times r}, \bar \C\in\Re^{r\times r\times r_3}$,
$V\in \Re^{n_2\times r_2}, W\in \Re^{n_3\times r_3}$ such that
\[\A=\bar \A\times_1 U,\quad \B=\bar \B \times_2 V,\quad \C=\bar\C \times_3 W. \]
Hence, $\X=\A\B\C=(\bar\A\bar\B\bar\C)\times_1 U\times_2 V\times_3 W$ according to (\ref{tripletucker}). From definition of the Tucker rank, we have the first inequalities of (\ref{eq3.16}-\ref{eq3.18}).

Assume that TriRank$(\A) = \bar r$.  Then there are tensors $\hat \A \in \Re^{n_1 \times \bar r \times \bar r}$, $\hat \B \in \Re^{\bar r \times r \times \bar r}$ and $\hat \C \in \Re^{\bar r \times \bar r \times r}$ such that $\A = \hat \A \hat \B \hat \C$.   Replacing $\X$ and $\A$ in the first inequality of
(\ref{eq3.16}) by $\A$ and $\hat A$, we have TucRank$(\A)_1 \le$  TucRank$(\hat \A)_1$.
Note that $\hat A \in \Re^{n_1 \times \bar r \times \bar r}$.  By the definition of the Tucker rank,
TucRank$(\hat \A)_1$ is the matrix rank of an $n_1 \times \bar r^2$ matrix.  Hence, TucRank$(\hat \A)_1 \le \bar r^2$.
This proves the second inequality of (\ref{eq3.16}).

Since $\A = \hat \A \hat \B \hat \C$ and $\A \in \Re^{n_1 \times r \times r}$, by (\ref{eq11}), TriRank$(\A) \le r =$ TriRank$(\X)$.   Then the third inequality of (\ref{eq3.16}) holds.

The second and third inequalities of (\ref{eq3.17}) and (\ref{eq3.18}) hold similarly.
\qed

\section{A Method for Checking The Triple Rank of a Third Order Tensor}

In this section, we present an algorithm for checking the triple rank of a third order tensor and establish its convergence.   Strictly speaking, our algorithm is not guaranteed to find the triple rank of a third order tensor $\X$ exactly. Instead, it gives an upper bound on the relative error obtainable by approximating $\X$ with a third order tensor $\A\B\C$ of triple rank not higher than a given integer $r$.    This algorithm will be useful in the next section to verify that  third order tensors from several practical datasets can be approximated by low triple rank tensors.


\subsection{A Modified Alternating Least Squares Method} \label{SubSect-ALSdecomp}


We are going to present a modified alternating least squares (MALS) algorithm for the triple decomposition of third order tensors in this subsection. Consider a given third order tensor $\X\in\Re^{n_1\times n_2\times n_3}$ with $n_1, n_2, n_3\ge1$ and a fixed positive integer $r \le$ mid$\{n_1, n_2, n_3\}$. The following cost function will be minimized
\begin{equation} \label{eq4.1}
  f(\A,\B,\C):=\left\|\X-\A\B\C\right\|_F^2=\sum_{i=1}^{n_1}\sum_{j=1}^{n_2}\sum_{t=1}^{n_3}\left( x_{ijt}-\sum_{p=1}^r\sum_{q=1}^r\sum_{s=1}^r a_{iqs}b_{pjs}c_{pqt} \right)^2,
\end{equation}
where $\A\in\Re^{n_1\times r\times r}$, $\B\in\Re^{r\times n_2\times r}$, $\C\in\Re^{r\times r\times n_3}$ are unknown.   In this way, we will obtain a triple decomposition $\A\B\C$ of triple rank not greater than
$r$, to approximate $\X$.

MALS is an iterative approach starting from an initial points $(\A^0,\B^0,\C^0)\in\Re^{n_1\times r\times r}\oplus\Re^{r\times n_2\times r} \oplus \Re^{r\times r\times n_3}$. We initialize $k\gets0$ and perform the following steps until the iterative sequence converges.

Update $\A^{k+1}$. Fixing $\B^k$ and $\C^k$, we solve a subproblem
\begin{equation*}
  \arg\min_{\A\in\Re^{n_1\times r\times r}}~\left\|\A\B^k\C^k-\X\right\|_F^2 + \lambda\left\|\A-\A^k\right\|_F^2,
\end{equation*}
where $\lambda >0$ is a constant in this algorithm. If $\lambda = 0$, then this is the classical ALS algorithm.  We take $\lambda > 0$.  Hence, our method is a modified ALS algorithm.  Let $A_{(1)}\in\Re^{n_1\times r^2}$ be the mode-1 unfolding of the tensor $\A$ and $X_{(1)}\in\Re^{n_1\times n_2n_3}$ be the mode-1 unfolding of the tensor $\X$. By introducing a matrix $F^k\in\Re^{r^2\times n_2n_3}$ with elements
\begin{equation}\label{Fk-def}
  F^k_{\ell m} = \sum_{p=1}^r b^k_{pjs}c^k_{pqt} \qquad
  \text{ where } \ell = q+(s-1)r, m = j+(t-1)n_2,
\end{equation}
the $\A$-subproblem may be represented as
\begin{eqnarray}\label{CloseF-A}
  \lefteqn{\arg\min_{A_{(1)}\in\Re^{n_1\times r^2}}~\left\|A_{(1)}F^k-X_{(1)}\right\|_F^2 + \lambda\left\|A_{(1)}-A_{(1)}^k\right\|_F^2} \nonumber\\
    &=& \left[X_{(1)}\left(F^k\right)^T+\lambda A_{(1)}^k\right]\left[F^k\left(F^k\right)^T+\lambda I_{r^2}\right]^{-1}.
\end{eqnarray}
Then, we obtain $\A^{k+1}$ from $A^{k+1}_{(1)}$ which is the closed-form solution \eqref{CloseF-A}.

Update $\B^{k+1}$. Consider the following subproblem
\begin{equation*}
  \arg\min_{\B\in\Re^{r\times n_2\times r}}~\left\|\A^{k+1}\B\C^k-\X\right\|_F^2 + \lambda\left\|\B-\B^k\right\|_F^2,
\end{equation*}
where $\A^{k+1}$ and $\C^k$ are known. Let $X_{(2)}\in\Re^{n_2\times n_1n_3}$ and $B_{(2)}\in\Re^{n_2\times r^2}$ be the $2$-mode unfolding of tensors $\X$ and $\B$, respectively. Define $G^k\in\Re^{r^2\times n_1n_3}$ with entries
\begin{equation}\label{Gk-def}
  G^k_{\ell m} = \sum_{q=1}^r a^{k+1}_{iqs}c^k_{pqt} \qquad\text{ where }
  \ell = p+(s-1)r, m = i+(t-1)n_1.
\end{equation}
Then, the $\B$-subproblem is rewritten as
\begin{eqnarray}\label{CloseF-B}
  \lefteqn{\arg\min_{B_{(2)}\in\Re^{n_2\times r^2}}~\left\|B_{(2)}G^k-X_{(2)}\right\|_F^2 +
  \lambda\left\|B_{(2)}-B_{(2)}^k\right\|_F^2} \nonumber\\
    &=& \left[X_{(2)}\left(G^k\right)^T+\lambda B_{(2)}^k\right]\left[G^k\left(G^k\right)^T+\lambda I_{r^2}\right]^{-1}.
\end{eqnarray}
Hence, $\B^{k+1}$ may be derived from $B^{k+1}_{(1)}$ defined by \eqref{CloseF-B}.

Update $\C^{k+1}$. Using $\A^{k+1}$ and $\B^{k+1}$ at hand, we minimize
\begin{equation*}
  \arg\min_{\C\in\Re^{r\times r\times n_3}}~\left\|\A^{k+1}\B^{k+1}\C-\X\right\|_F^2 + \lambda\left\|\C-\C^k\right\|_F^2.
\end{equation*}
Let $H\in\Re^{r^2\times n_1n_2}$ be a matrix with entries
\begin{equation}\label{Hk-def}
  H^k_{\ell m} = \sum_{s=1}^r a^{k+1}_{iqs}b^{k+1}_{pjs} \qquad\text{ where }
  \ell = p+(q-1)r, m = i+(j-1)n_1.
\end{equation}
Then, we derive
\begin{eqnarray}\label{CloseF-C}
  \lefteqn{\arg\min_{C_{(3)}\in\Re^{n_3\times r^2}}~\left\|C_{(3)}H^k-X_{(3)}\right\|_F^2 + \lambda\left\|C_{(3)}-C_{(3)}^k\right\|_F^2} \nonumber\\
    &=& \left[X_{(3)}\left(H^k\right)^T+\lambda C_{(3)}^k\right]\left[H^k\left(H^k\right)^T+\lambda I_{r^2}\right]^{-1},
\end{eqnarray}
where $X_{(3)}$ and $C_{(3)}$ are the 3-mode unfolding of $\X$ and $\C$, respectively.
The third order tensor $\C^{k+1}$ is a tensor-form of \eqref{CloseF-C}.

\begin{algorithm}
\caption{Modified Alternating Least Squares (MALS) algorithm for triple decomposition.}\label{ALSa}
\begin{algorithmic}[1]
  \STATE Set $\gamma\in[1,2)$ and $\lambda>0$. Choose an integer $r\ge1$ and an initial point $\A^0\in\Re^{n_1\times r\times r}$,
    $\B^0\in\Re^{r\times n_2\times r}$, and $\C^0\in\Re^{r\times r\times n_3}$. Set $k\gets 0$.
  \STATE Compute $\widetilde{\A}^k$ by \eqref{CloseF-A} and set $\A^{k+1}=\gamma\widetilde{\A}^k+(1-\gamma)\A^k$.
  \STATE Compute $\widetilde{\B}^k$ by \eqref{CloseF-B} and set $\B^{k+1}=\gamma\widetilde{\B}^k+(1-\gamma)\B^k$.
  \STATE Compute $\widetilde{\C}^k$ by \eqref{CloseF-C} and set $\C^{k+1}=\gamma\widetilde{\C}^k+(1-\gamma)\C^k$.
  \STATE Set $k\gets k+1$ and goto Step 2.
\end{algorithmic}
\end{algorithm}

Set $k\gets k+1$ and repeat the process. See Algorithm \ref{ALSa} for a complete algorithm. Here, we use the extrapolation technique with step size $\gamma\in[1,2)$ to deal with the swamp effect. ALS may terminate if the difference between two iterates is small enough, i.e.,
\begin{equation*}
  \max\left\{\frac{\|\A^{k+1}-\A^{k}\|_F}{\|\A^{k+1}\|_F}, \frac{\|\B^{k+1}-\B^{k}\|_F}{\|\B^{k+1}\|_F}, \frac{\|\C^{k+1}-\C^{k}\|_F}{\|\C^{k+1}\|_F}\right\} \le \varepsilon,
\end{equation*}
or the iteration arrives a preset maximal iterative number.

\subsection{Convergence Analysis of Algorithm \ref{ALSa}}

If $\gamma = 1$ in Algorithm \ref{ALSa}, then this algorithm can be regarded as a special case of the block coordinate descent (BCD) method with proximal update, studied in \cite{XY13}.   However, $\gamma > 1$ is an extrapolation step, which will speed the algorithm.   In the convergence analysis of \cite{XY13},
the Kurdyka-{-L}ojasiewicz (KL) inequality is assumed to hold at a limiting point of the iterative sequence.  We do not need to make this assumption.   Hence, the convergence result of \cite{XY13} cannot cover Algorithm \ref{ALSa}.

Our algorithm is closer to the seminorm regularized alternating least squares (SRALS) algorithm for CP tensor decomposition presented in \cite{CSXY19}.  There are two differences between our algorithm and the SRALS algorithm.  First, SRALS is for CP decomposition, while Algorithm \ref{ALSa} is for triple decomposition.  Second, Step 8 of SRALS is not used in Algorithm \ref{ALSa}, as it is not necessary here.   Otherwise,
Algorithm \ref{ALSa} is similar to SRALS.

Note that as we set $\lambda > 0$, our algorithm is not an ALS method, but a modified ALS method.  An argument following Lemma 4.2 of \cite{CSXY19} justifies this.

Then, with a similar argument for the proofs of Theorems 4.1 and 4.4, Lemma 4.8 and Theorem 4.9 of \cite{CSXY19}, we have Theorem \ref{t4.1}.  Before stating the theorem, we first state the definition of the well-known KL inequality \cite{CSXY19}.

\begin{definition}
Let $f: U \to \Re$, where $U \subseteq \Re^N$ is an open set, be an analytic function.    Let $\x \in U$.   We say that the Kurdyka-{-L}ojasiewicz (KL) inequality holds at $\x$ if there is a neighborhood $V$ of $\x$, an exponent $\theta \in [{1 \over 2}, 1)$ and a constant $C$ such that for any $\y \in V$,
$$|f(\y) - f(\x)|^\theta \le C\|\nabla f(\y)\|_2.$$
\end{definition}

Then we have the following theorem as stated early.

\begin{theorem} \label{t4.1}
Denote $\X^k = \A^k\B^k\C^k$ in Algorithm \ref{ALSa}.  Let $f$ be defined in (\ref{eq4.1}).   Suppose that
Algorithm \ref{ALSa} generates a sequence $\{ \X^k \}$.   If $\X^k = \X^{k+1}$ for some $k$, then $\X^k$ is a critical point of $f$.  Otherwise, an infinite sequence $\{ \X^k \}$ is generated.   If this sequence is bounded, then this sequence converges to a critical point $\bar \X$ of $f$, and the KL inequality holds at $\bar \X$.  If $\theta = {1 \over 2}$ in the KL inequality, then there exists $c > 0$ and $Q \in [0, 1)$ such that
$$\| \X^k - \bar \X \|_F \le c Q^k.$$
If $\theta = ({1 \over 2}, 1)$ in the KL inequality, then there exists $c > 0$ such that
$$\| \X^k - \bar \X \|_F \le c k^{-{1-\theta \over 2\theta - 1}}.$$
\end{theorem}


\section{Practical Data of Third Order Tensors}


In this section, we investigate practical data from applications and show that they can be approximated by triple decomposition of low triple ranks very well.

\subsection{Abilene Internet Traffic Data}

The first application we consider is the internet traffic data.   The data set is the Abilene data set \footnote{The abilene observatory data collections. http://abilene.internet2.edu/observatory/data-collections.html} \cite{XWWXWZ16}.

\begin{figure}
  \centering
  \includegraphics[width=0.8\textwidth]{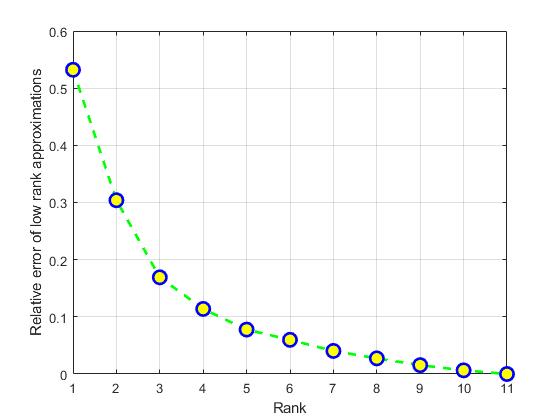}\\
  \caption{Relative errors of low triple rank approximations of the $11\times 11\times 2016$ internet traffic tensor from Abilene dataset.}\label{Abin-1}
\end{figure}

The Abilene data arises from the backbone network located in North America. There are 11 routers: Atlanta GA, Chicago IL, Denver CO, Houston TX, Indianapolis, Kansas City MO, Los Angeles CA, New York NY, Sunnyvale CA, Seattle WA, and Washington DC. These routers send and receive data. Thus we get 121 original-destination (OD) pairs. For each OD pair, we record the internet traffic data of every 5 minutes in a week from Dec 8, 2003 to Dec 14, 2003. Hence, there are $7\times 24 \times 60/5 = 2016$ numbers for each OD pairs. In this way, we get a third order tensor $\X_{Abil}$ with size $11$-by-$11$-by-$2016$.  This model was used in \cite{ADKM11, ZZXC15} for internet traffic data recovery.

Now, we examine the triple decomposition approximation of the tensor $\X_{Abil}\in\Re^{11\times 11\times 2016}$ with different triple rank upper bound among $1$ to $11$. For each triple rank upper bound, we compute the triple decomposition approximation $\A\B\C$ by Algorithm \ref{ALSa} and calculate the relative error of low triple rank approximation
\begin{equation*}
  \mathrm{Relative Error} = \frac{\|\X_{Abil}-\A\B\C\|_F}{\|\X_{Abil}\|_F}.
\end{equation*}
Figure \ref{Abin-1} illustrates the relative error of the low rank approximations via triple rank upper bound $r$. When we take $r=5$ and $r=7$, the relative error is about $7.8\%$ and $4.0\%$, respectively.    This
shows that the Abilene data can ba approximated by triple decomposition of low triple rank well.
Obviously, the relative error is zero if $r=11$ as this is an upper bound on the rank as shown in \eqref{eq11}.

\begin{figure}
  \centering
  \includegraphics[width=0.8\textwidth]{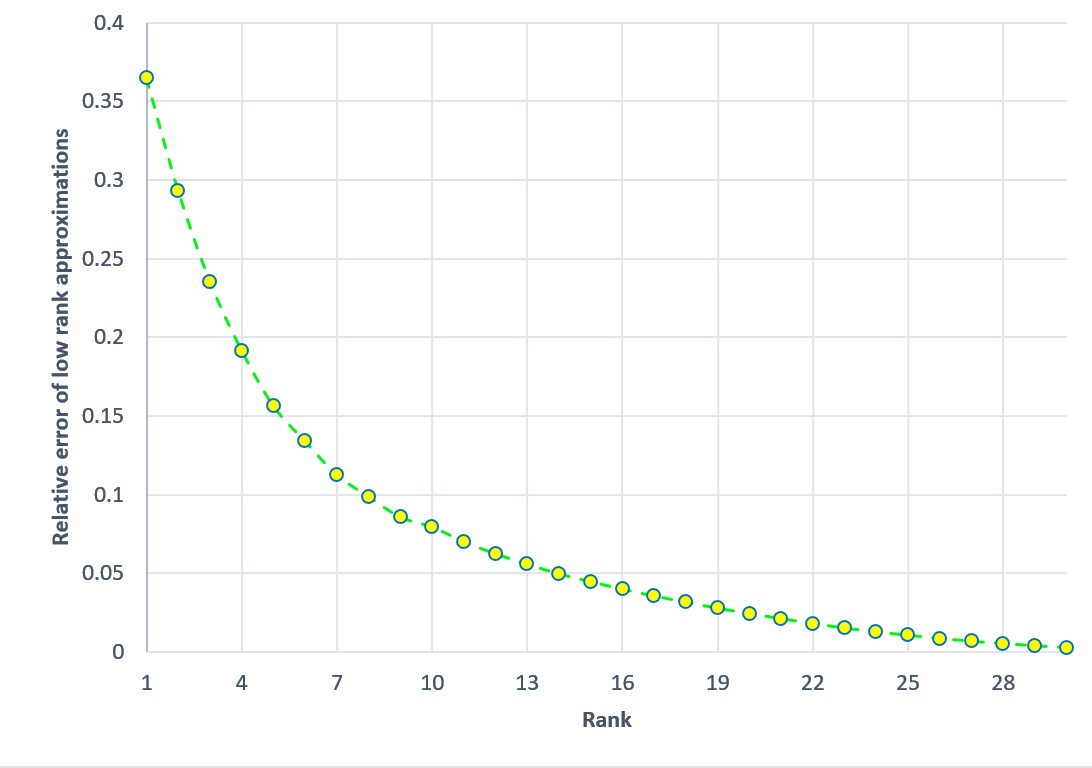}\\
  \caption{Relative errors of low triple rank approximations of the $121\times 96\times 21$ internet traffic tensor from Abilene dataset.}\label{Abin-2}
\end{figure}

A similar conclusion is obtained if we view the Abilene traffic data as a third order tensor arranged differently as $\tilde{\X}_{Abil}\in\Re^{121\times 96\times 21}$, which is indexed by $121$ source-destination pairs, $96$ time slots for each day, and $21$ days.   This is the model used in \cite{XWWXWZ16}.   Figure~\ref{Abin-2} shows the relative error of the low rank approximations obtained by Algorithm \ref{ALSa} as a function of the target ranks upto  30.   Actually, this is more illustrative as here $n_2 = 96$ and the low rank triple decomposition is very good when $r \ge 25$.

\subsection{ORL Face Data}

We now investigate the ORL face data in AT \& T Laboratories Cambridge \cite{CSZZC19, SH94, XY13}.

\begin{figure}
  \centering
  \includegraphics[width=0.8\textwidth]{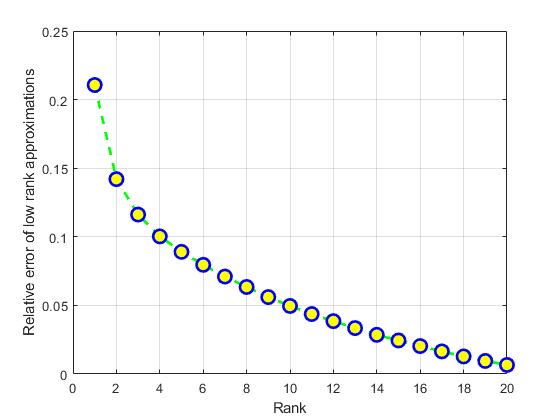}\\
  \caption{Low rank approximations of a third order ORL face data tensor of size $112\times 92\times 10$.}\label{ORL-A}
\end{figure}

\begin{figure}
  \centering
  \includegraphics[width=\textwidth]{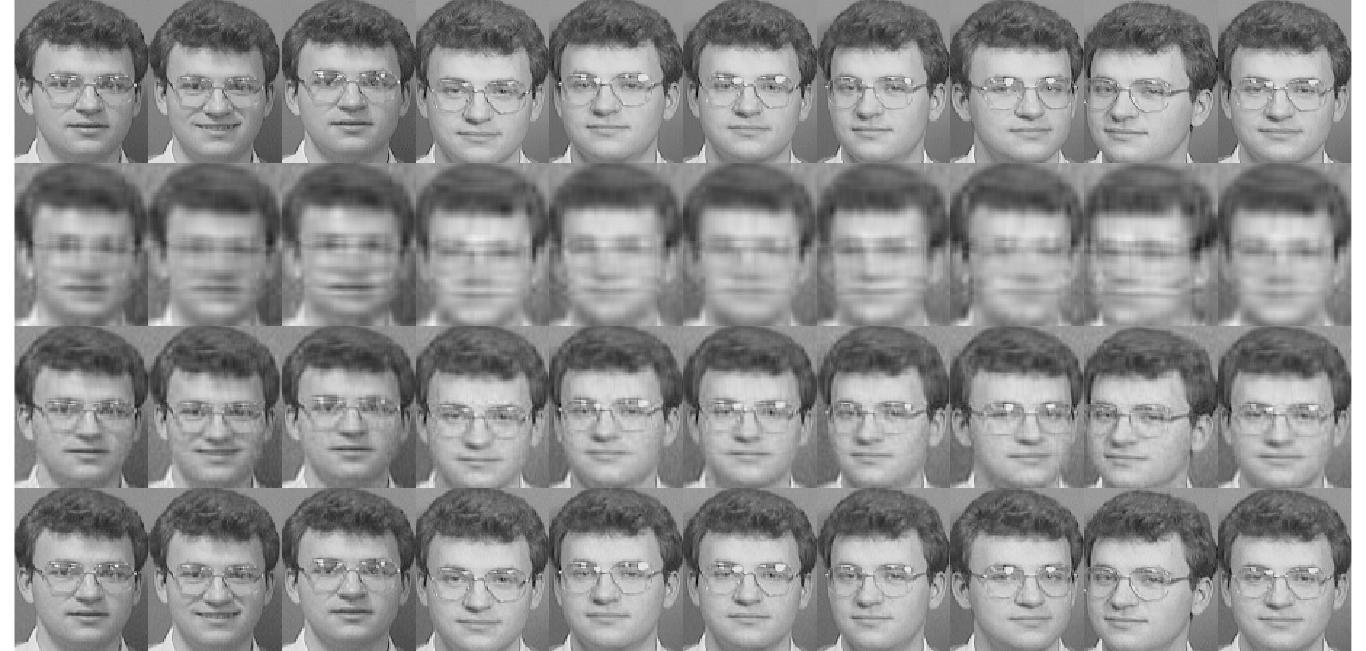}\\
  \caption{Illustration of faces from the ORL dataset. Original images are illustrated in the first line. Approximations with rank 4, 10, 16 are shown in lines two, three, four, respectively. }\label{ORL_B}
\end{figure}

The ORL dataset of faces contains images of 40 persons. Each image has $112\times 92$ pixels. For each person, there are 10 images taken at different times, varying the lighting, facial expressions and facial details. For instance, the first line of Figure \ref{ORL_B} illustrates ten images of a person. Hence, there is a $112$-by-$92$-by-$10$ tensor $\T_{face}$. Using Algorithm \ref{ALSa}, we compute best low triple rank approximations of the tensor $\T_{face}$. The relative error of approximations via triple ranks are illustrated in Figure \ref{ORL-A}. When the triple rank upper bound $r=4,10,16$, the relative error of low triple rank approximations are $10.03\%,4.96\%,2.04\%$, respectively. Corresponding images of low triple rank approximations are illustrated in lines 2--4 of Figure \ref{ORL_B}.

This result shows clearly the ORL data can be approximated by low rank triple decomposition very well.

\section{A Tensor Recovery Method and Its Convergence Analysis}


In this section, we consider the tensor recovery problem:
\begin{equation}\label{Ten-Recovery}
  \min~ \|\mathbb{P}(\A\B\C)-\dd\|^2_F,
\end{equation}
where $\mathbb{P}$ is a linear operator, $\dd\in\Re^m$ is a given vector, and $\A\in\Re^{n_1\times r\times r}$, $\B\in\Re^{r\times n_2\times r}$, and $\C\in\Re^{r\times r\times n_3}$ are unknown. To solve \eqref{Ten-Recovery}, we introduce a surrogate tensor $\X = (x_{ijt}) \in\Re^{n_1\times n_2\times n_3}$ and transform \eqref{Ten-Recovery} to the closely related optimization problem
\begin{equation}\label{Ten-Recov-1}
\begin{aligned}
  \min~ & f(\X,\A,\B,\C):=\|\A\B\C-\X\|_F^2, \\
  \st~ & \mathbb{P}(\X)=\dd.
\end{aligned}
\end{equation}
Here, we slightly abuse the notation $f$ for denoting an objective function.

\subsection{A Tensor Recovery Method}

We propose a modified alternating least squares algorithm for solving the tensor recovery problem \eqref{Ten-Recov-1}. For a fixed positive integer $r$, we choose $\A^0\in\Re^{n_1\times r\times r}$, $\B^0\in\Re^{r\times n_2\times r}$, $\C^0\in\Re^{r\times r\times n_3}$, $\X^0\in\Re^{n_1\times n_2\times n_3}$ and set $k\gets0$.
Using an approach similar to that introduced in subsection \ref{SubSect-ALSdecomp}, we perform the following steps.

Update $\X^{k+1}$. We solve a subproblem
\begin{equation*}
\begin{aligned}
  \arg\min_{\X\in\Re^{n_1\times n_2\times n_3}}~& \|\X-\A^k\B^k\C^k\|_F^2+\lambda\|\X-\X^k\|_F^2 \\
  \st~ & \mathbb{P}(\X)=\dd.
\end{aligned}
\end{equation*}
That is
\begin{equation*}
\begin{aligned}
  \arg\min_{\X\in\Re^{n_1\times n_2\times n_3}}~& \left\|\X-\tfrac{1}{1+\lambda}\left(\A^k\B^k\C^k+\lambda\X^k\right)\right\|_F^2 \\
  \st~ & \mathbb{P}(\X)=\dd.
\end{aligned}
\end{equation*}
Define an operator $\mathrm{vec}:\Re^{n_1\times n_2\times n_3}\to\Re^{n_1n_2n_3}$ that maps $x_{ijt}$ to $\hat x_{\ell}$ where $\ell = i+(j-1)n_1+(t-1)n_1n_2$. Then, the equality constraint $\mathbb{P}(\X)=\dd$ may be rewritten as $P\mathrm{vec}(\X)=\dd$, where $P$ is the $m$-by-$(n_1n_2n_3)$ matrix corresponding to the application of the operator $\mathbb{P}$ when viewed as a linear transformation from $\mathrm{vec}(\X)$ to $\dd$.   Here we assume that $PP^T$ is invertible.  Thus, the above optimization problem may be represented as
\begin{equation*}
\begin{aligned}
  \arg\min~& \left\|\mathrm{vec}(\X)-\tfrac{1}{1+\lambda}\mathrm{vec}\left(\A^k\B^k\C^k+\lambda\X^k\right)\right\|_F^2 \\
  \st~ & P\mathrm{vec}(\X)=\dd,
\end{aligned}
\end{equation*}
which has a closed-form solution
\begin{equation}\label{CloseF-T2}
  \left[I-P^T\left(PP^T\right)^{-1}P\right]\tfrac{1}{1+\lambda}\mathrm{vec}\left(\A^k\B^k\C^k+\lambda\X^k\right)+P^T\left(PP^T\right)^{-1}\dd.
\end{equation}
This is defined as $\mathrm{vec}(\widetilde{\X}^k)$. Next, we set
\begin{equation*}
  \X^{k+1} = \gamma\widetilde{\X}^k + (1-\gamma)\X^k.
\end{equation*}

Update $\A^{k+1}$. To solve
\begin{equation*}
  \arg\min_{\A\in\Re^{n_1\times r\times r}}~\|\A\B^k\C^k-\X^{k+1}\|_F^2 + \lambda\|\A-\A^k\|_F^2,
\end{equation*}
we obtain $\widetilde{\A}^k$ by calculating $\widetilde{A}^k_{(1)}$ which is
\begin{eqnarray}\label{CloseF-A2}
  \lefteqn{\arg\min_{A_{(1)}\in\Re^{n_1\times r^2}}~\|A_{(1)}F^k-X^{k+1}_{(1)}\|_F^2 + \lambda\|A_{(1)}-A_{(1)}^k\|_F^2} \nonumber\\
    &=& \left[X^{k+1}_{(1)}\left(F^k\right)^T+\lambda A_{(1)}^k\right]\left[F^k\left(F^k\right)^T+\lambda I_{r^2}\right]^{-1},
\end{eqnarray}
where $F^k$ is defined in \eqref{Fk-def} using $\B^k$ and $\C^k$, and $A_{(1)}$, $A^k_{(1)}$, $X^{k+1}_{(1)}$ are 1-mode unfolding of tensors $\A$, $\A^k$, $\X^{k+1}$, respectively.   We apply extrapolation to set
\begin{equation*}
  \A^{k+1} = \gamma\widetilde{\A}^k + (1-\gamma)\A^k.
\end{equation*}

Update $\B^{k+1}$. To solve
\begin{equation*}
  \arg\min_{\B\in\Re^{r\times n_2\times r}}~\left\|\A^{k+1}\B\C^k-\X^{k+1}\right\|_F^2 + \lambda\left\|\B-\B^k\right\|_F^2,
\end{equation*}
we obtain $\widetilde{\B}^k$ by calculating $\widetilde{B}^k_{(2)}$ which is
\begin{eqnarray}\label{CloseF-B2}
  \lefteqn{\arg\min_{B_{(2)}\in\Re^{n_2\times r^2}}~\left\|B_{(2)}G^k-X^{k+1}_{(2)}\right\|_F^2 + \lambda\left\|B_{(2)}-B_{(2)}^k\right\|_F^2} \nonumber\\
    &=& \left[X^{k+1}_{(2)}\left(G^k\right)^T+\lambda B_{(2)}^k\right]\left[G^k\left(G^k\right)^T+\lambda I_{r^2}\right]^{-1},
\end{eqnarray}
where $G^k$ is defined in \eqref{Gk-def} using $\A^{k+1}$ and $\C^k$, and $B_{(2)}$, $B^k_{(2)}$, $X^{k+1}_{(2)}$ are $2$-mode unfolding of tensors $\B$, $\B^k$, $\X^{k+1}$, respectively. We apply extrapolation to set
\begin{equation*}
  \B^{k+1} = \gamma\widetilde{\B}^k + (1-\gamma)\B^k.
\end{equation*}

Update $\C^{k+1}$. To solve
\begin{equation*}
  \arg\min_{\C\in\Re^{r\times r\times n_3}}~\|\A^{k+1}\B^{k+1}\C-\X^{k+1}\|_F^2 + \lambda\|\C-\C^k\|_F^2,
\end{equation*}
we obtain $\widetilde{\C}^k$ by calculating $\widetilde{C}^k_{(3)}$ which is
\begin{eqnarray}\label{CloseF-C2}
  \lefteqn{\arg\min_{C_{(3)}\in\Re^{n_3\times r^2}}~\|C_{(3)}H^k-X^{k+1}_{(3)}\|_F^2 + \lambda\|C_{(3)}-C_{(3)}^k\|_F^2} \nonumber\\
    &=& \left[X^{k+1}_{(3)}\left(H^k\right)^T +\lambda C_{(3)}^k\right]\left[H^k\left(H^k\right)^T+\lambda I_{r^2}\right]^{-1},
\end{eqnarray}
where $H^k$ is defined in \eqref{Hk-def} using $\A^{k+1}$ and $\B^{k+1}$, and $C_{(3)}$, $C^k_{(3)}$, $X^{k+1}_{(3)}$ are 3-mode unfolding of tensors $\C$, $\C^k$, $\X^{k+1}$, respectively.  We apply extrapolation to set
\begin{equation*}
  \C^{k+1} = \gamma\widetilde{\C}^k + (1-\gamma)\C^k.
\end{equation*}

Subsequently, we set $k\gets k+1$ and repeat this process until convergence. The detailed algorithm are illustrated in Algorithm \ref{ALSb}.

\begin{algorithm}
\caption{Modified Alternating Least Squares (MALS) algorithm for recovering a third order tensor.}\label{ALSb}
\begin{algorithmic}[1]
  \STATE Set $\gamma\in[1,2)$ and $\lambda>0$. Choose an integer $r\ge1$ and an initial point $\A^0\in\Re^{n_1\times r\times r}$,
    $\B^0\in\Re^{r\times n_2\times r}$, $\C^0\in\Re^{r\times r\times n_3}$, and $\X^0\in\Re^{n_1\times n_2\times n_3}$. Set $k\gets 0$.
  \STATE Compute $\widetilde{\X}^k$ by \eqref{CloseF-T2} and set $\X^{k+1}=\gamma\widetilde{\X}^k+(1-\gamma)\X^k$.
  \STATE Compute $\widetilde{\A}^k$ by \eqref{CloseF-A2} and set $\A^{k+1}=\gamma\widetilde{\A}^k+(1-\gamma)\A^k$.
  \STATE Compute $\widetilde{\B}^k$ by \eqref{CloseF-B2} and set $\B^{k+1}=\gamma\widetilde{\B}^k+(1-\gamma)\B^k$.
  \STATE Compute $\widetilde{\C}^k$ by \eqref{CloseF-C2} and set $\C^{k+1}=\gamma\widetilde{\C}^k+(1-\gamma)\C^k$.
  \STATE Set $k\gets k+1$ and goto Step 2.
\end{algorithmic}
\end{algorithm}

\subsection{Convergence Analysis}
\newcommand{\vt}[1]{\mathrm{vec}(#1)}

We now present convergence analysis for this algorithm. For convenience, we collect all variables as a undetermined vector
\begin{equation*}
  \y := \left(\vt{\X}^T,\vt{A_{(1)}}^T,\vt{B_{(2)}}^T,\vt{C_{(3)}}^T\right)^T \in\Re^{n_1n_2n_3+(n_1+n_2+n_3)r^2}.
\end{equation*}
The feasible region of $\y$ is defined by
\begin{equation*}
  \Omega:=\{\vt{\X}\in\Re^{n_1n_2n_3}:P\vt{\X}=\dd\}\oplus\Re^{n_1r^2}\oplus\Re^{n_2r^2}\oplus\Re^{n_3r^2}.
\end{equation*}
We analyze the convergence of Algorithm \ref{ALSb} the solving an optimization problem
\begin{equation}\label{CA-1}
  \min~f(\y):=f(\X,\A,\B,\C) = \|\X-\A\B\C\|_F^2 \qquad \mathrm{s.t.}~\y\in\Omega.
\end{equation}
To simplify notation, we use $\y=(\X,\A,\B,\C)$ in the following analysis.

By optimization theory, $\y^*$ is a stationary point of \eqref{CA-1} if and only if the projected negative gradient of $f$ at $\y^*$ vanishes. In the following, we derive the formula of the projected gradient of $f$. First, let $\y=(\X,\A,\B,\C)\in\Omega$. Since
\begin{eqnarray*}
  f(\y) &=& \|\vt{\X}-\vt{\A\B\C}\|^2  \\
    &=& \langle\vt{\X},\vt{\X}\rangle-2\langle\vt{\X},\vt{\A\B\C}\rangle+\langle\vt{\A\B\C},\vt{\A\B\C}\rangle.
\end{eqnarray*}
Hence, $\nabla_{\vt{\X}}f = 2\vt{\X}-2\vt{\A\B\C} = 2\vt{\X-\A\B\C}$. Since the set $\{\vt{\X}\in\Re^{n_1n_2n_3}:P\vt{\X}=\dd\}$ is an affine manifold, we obtain the projected gradient of $\X$-part
\begin{equation*}
  \left[I-P^T\left(PP^T\right)^{-1}P\right](2\vt{\X-\A\B\C})= 2 \left[I-P^T\left(PP^T\right)^{-1}P\right]\vt{\X-\A\B\C}
\end{equation*}
directly.

Next, we rewrite $f(\y)$ as
\begin{eqnarray*}
  f(\y) &=& \|A_{(1)}F-X_{(1)}\|_F^2 \\
    &=& \langle A_{(1)}F,A_{(1)}F \rangle -2\langle A_{(1)}F,X_{(1)} \rangle +\langle X_{(1)},X_{(1)} \rangle \\
    &=& \langle A_{(1)},A_{(1)}FF^T \rangle -2\langle A_{(1)},X_{(1)}F^T \rangle +\langle X_{(1)},X_{(1)} \rangle,
\end{eqnarray*}
where $F$ is defined by \eqref{Fk-def} using $\B$ and $\C$. Hence, the $A_{(1)}$-part of the (projected) gradient is
\begin{equation*}
  2\left(A_{(1)}FF^T-X_{(1)}F^T\right)  = 2(A_{(1)}F-X_{(1)})F^T.
\end{equation*}
We may write the $\vt{A_{(1)}}$-part of the gradient in $2\vt{(A_{(1)}F-X_{(1)})F^T}$ to corresponding to the vector form on $\vt{A_{(1)}}$.

Finally, by a similar approach, the $B_{(2)}$-part and the $C_{(3)}$-part of the (projected) gradient are
\begin{equation*}
  2(B_{(2)}G-X_{(2)})G^T, \qquad 2(C_{(3)}H-X_{(3)})H^T,
\end{equation*}
respectively. Here, $G$ is defined by \eqref{Gk-def} using $\A$ and $\C$; and $H$ is defined by \eqref{Hk-def} using $\A$ and $\B$.
Therefore, we get the projected gradient of $f$ at $\y=(\X,\A,\B,\C)\in\Omega$:
\begin{equation}\label{pgy}
  \Pi_{\Omega}\left(\nabla f(\y)\right) = 2\left(\begin{array}{c}
    [I-P^T\left(PP^T\right)^{-1}P]\vt{\X-\A\B\C} \\
    \vt{(A_{(1)}F-X_{(1)})F^T} \\
    \vt{(B_{(2)}G-X_{(2)})G^T} \\
    \vt{(C_{(3)}H-X_{(3)})H^T} \\
  \end{array}\right),
\end{equation}
where $\Pi_{\Omega}(\cdot)$ denotes the projection onto the feasible $\Omega$.
We have the following lemma on the optimality condition.

\begin{lemma}\label{Lem-5.1}
  Let $\y^*=  \left(\X^*,\A^*,\B^*,\C^*\right)^T \in\Omega$ be the optimal solution of optimization problem \eqref{CA-1}. Then, the projected negative gradient of $f$ at $\y^*$ vanishes, i.e.,
  \begin{eqnarray}
    & (I-P^T(PP^T)^{-1}P)\mathrm{vec}(\A^*\B^*\C^*-\X^*)=0, & \\
    & (X_{(1)}^*-A_{(1)}^*F^*)(F^*)^T = 0, & \\
    & (X_{(2)}^*-B_{(2)}^*G^*)(G^*)^T = 0, & \\
    & (X_{(3)}^*-C_{(3)}^*H^*)(H^*)^T = 0, &
  \end{eqnarray}
  where $F^*$ is defined by \eqref{Fk-def} using $\B^*$ and $\C^*$, $G^*$ is defined by \eqref{Gk-def} using $\A^*$ and $\C^*$, $H^*$ is defined by \eqref{Hk-def} using $\A^*$ and $\B^*$. That is to say, $\y^*=(\X^*,\A^*,\B^*,\C^*)$ is a stationary point of \eqref{CA-1}.
\end{lemma}

Now, we consider the case that the sequence generated by Algorithm \ref{ALSb} converges in a finite number of iterations.

\begin{lemma}\label{Lem-5.2}
  If there exists an iteration $k$ such that $\y^k=\y^{k+1}$, i.e.,
  \begin{equation*}
    (\X^k,\A^k,\B^k,\C^k)=(\X^{k+1},\A^{k+1},\B^{k+1},\C^{k+1}),
  \end{equation*}
  then $(\X^k,\A^k,\B^k,\C^k)$ is a stationary point of \eqref{CA-1}.
\end{lemma}
{\bf Proof}
  First, for the $\X$-part, since $\widetilde{\X}^k$ is generated by \eqref{CloseF-T2}, we know
  \begin{equation*}
    \mathrm{vec}(\widetilde{\X}^k) = \left(I-P^T(PP^T)^{-1}P\right)\tfrac{1}{1+\lambda}\mathrm{vec}(\A^k\B^k\C^k+\lambda\X^k) + P^T(PP^T)^{-1}\dd.
  \end{equation*}
  In addition, because $\X^k$ satisfies $P\vt{\X^k}=\dd$, it yields that
  \begin{equation*}
    \mathrm{vec}(\X^k) = \left(I-P^T(PP^T)^{-1}P\right)\mathrm{vec}(\X^k) + P^T(PP^T)^{-1}\dd.
  \end{equation*}
  Combining the above two equations, we have
  \begin{eqnarray}
    \mathrm{vec}(\widetilde{\X}^k-\X^k) &=& \left(I-P^T(PP^T)^{-1}P\right)\tfrac{1}{1+\lambda}\mathrm{vec}(\A^k\B^k\C^k-\X^k) \nonumber\\
      &=& \tfrac{1}{1+\lambda}\left(I-P^T(PP^T)^{-1}P\right)\mathrm{vec}(\A^k\B^k\C^k-\X^k). \label{xxx}
  \end{eqnarray}
  From $\X^k=\X^{k+1}$, we get $\widetilde{\X}^k-\X^k=0$. Hence, $\left(I-P^T(PP^T)^{-1}P\right)\mathrm{vec}(\A^k\B^k\C^k-\X^k)=0$, i.e., the $\X$-part of the projected negative gradient of $f$ vanishes.

  It yields from \eqref{CloseF-A2} that
  \begin{equation*}
    \widetilde{A}_{(1)}^k = \left(X^{k+1}_{(1)}\left(F^k\right)^T+\lambda A_{(1)}^k\right)\left(F^k\left(F^k\right)^T+\lambda I_{r^2}\right)^{-1},
  \end{equation*}
  which implies
  \begin{equation}\label{aaa}
    (X_{(1)}^{k+1}-A_{(1)}^kF^k)\left(F^k\right)^T = (\widetilde{A}_{(1)}^k-A_{(1)}^k)\left(F^k\left(F^k\right)^T+\lambda I_{r^2}\right).
  \end{equation}
  By $\A^k=\A^{k+1}$, we have $\widetilde{\A}^k=\A^k$. Hence $(X^k_{(1)}-A_{(1)}^kF^k)(F^k)^T=0$, i.e., the $\A$-part of the projected negative gradient of $f$ vanishes.

  Finally, by a similar discussion, we find that the $\B$-part and the $\C$-part of the projected negative gradient of $f$ also vanish.
  Hence, by Lemma \ref{Lem-5.1}, $\y^k=(\X^k,\A^k,\B^k,\C^k)$ is a stationary point of \eqref{CA-1}.
\qed

Next, we consider the case that Algorithm \ref{ALSb} generates an infinite sequence of iterates.

\begin{lemma}\label{Lem-5.3}
  Let $\{\y^k\}_{k=0,1,2,\dots}$ be a sequence of iterates generated by Algorithm \ref{ALSb}. Then, we have
  \begin{eqnarray*}
    f(\y^k) - f(\y^{k+1}) \ge \frac{2\lambda}{\gamma}\|\y^k-\y^{k+1}\|^2.
  \end{eqnarray*}
\end{lemma}
{\bf Proof} For the $\X$-part, we have
  \begin{eqnarray*}
    \lefteqn{ f(\X^k,\A^k,\B^k,\C^k)-f(\X^{k+1},\A^k,\B^k,\C^k) } \\
      &=& \|\X^k-\A^k\B^k\C^k\|_F^2 - \|\gamma\widetilde{\X}^k+(1-\gamma)\X^k-\A^k\B^k\C^k\|_F^2 \\
      &=& 2\langle \A^k\B^k\C^k-\X^k, \gamma(\widetilde{\X}^k-\X^k) \rangle - \|\gamma(\widetilde{\X}^k-\X^k)\|_F^2 \\
      &=& 2\gamma\langle \mathrm{vec}(\A^k\B^k\C^k-\X^k), \mathrm{vec}(\widetilde{\X}^k-\X^k) \rangle - \gamma^2\|\widetilde{\X}^k-\X^k\|_F^2 \\
      &=& 2\gamma\langle \mathrm{vec}(\A^k\B^k\C^k-\X^k), \tfrac{1}{1+\lambda}\left(I-P^T(PP^T)^{-1}P\right)\mathrm{vec}(\A^k\B^k\C^k-\X^k) \rangle - \gamma^2\|\widetilde{\X}^k-\X^k\|_F^2 \\
      &=& 2\gamma(1+\lambda)\|\tfrac{1}{1+\lambda}\left(I-P^T(PP^T)^{-1}P\right)\mathrm{vec}(\A^k\B^k\C^k-\X^k)\|^2 - \gamma^2\|\widetilde{\X}^k-\X^k\|_F^2 \\
      &=& 2\gamma(1+\lambda)\|\mathrm{vec}(\widetilde{\X}^k-\X^k)\|^2 - \gamma^2\|\widetilde{\X}^k-\X^k\|_F^2 \\
      &=& (2\gamma(1+\lambda)-\gamma^2)\|\widetilde{\X}^k-\X^k\|_F^2 \\
      &\ge& 2\lambda\gamma\|\widetilde{\X}^k-\X^k\|_F^2 \\
      &=& \frac{2\lambda}{\gamma}\|\X^{k+1}-\X^k\|_F^2,
  \end{eqnarray*}
  where the fourth and the sixth equalities hold by \eqref{xxx}, the fifth equality holds because $I-P^T(PP^T)^{-1}P$ is idempotent matrix, i.e., $(I-P^T(PP^T)^{-1}P)^2=I-P^T(PP^T)^{-1}P$, and the last inequality holds since $2\gamma>\gamma^2$.

  For the $\A$-part, we could establish
  \begin{eqnarray*}
    \lefteqn{ f(\X^{k+1},\A^k,\B^k,\C^k)-f(\X^{k+1},\A^{k+1},\B^k,\C^k) } \\
      &=& \|\A^k\B^k\C^k-\X^{k+1}\|_F^2 - \|(\gamma\widetilde{\A}^k+(1-\gamma)\A^k)\B^k\C^k-\X^{k+1}\|_F^2 \\
      &=& -2\langle \gamma(\widetilde{\A}^k-\A^k)\B^k\C^k, \A^k\B^k\C^k-\X^{k+1} \rangle -\|\gamma(\widetilde{\A}^k-\A^k)\B^k\C^k\|_F^2 \\
      &=& 2\gamma\langle (\widetilde{A}_{(1)}^k-A_{(1)}^k)F^k, X_{(1)}^{k+1}-A_{(1)}^kF^k \rangle -\gamma^2\|(\widetilde{A}_{(1)}^k-A_{(1)}^k)F^k\|_F^2 \\
      &=& 2\gamma\langle (\widetilde{A}_{(1)}^k-A_{(1)}^k), (X_{(1)}^{k+1}-A_{(1)}^kF^k){F^k}^T \rangle -\gamma^2\|(\widetilde{A}_{(1)}^k-A_{(1)}^k)F^k\|_F^2 \\
      &=& 2\gamma\langle (\widetilde{A}_{(1)}^k-A_{(1)}^k), (\widetilde{A}_{(1)}^k-A_{(1)}^k)(F^k(F^k)^T+\lambda I_{r^2}) \rangle -\gamma^2\|(\widetilde{A}_{(1)}^k-A_{(1)}^k)F^k\|_F^2 \\
      &=& 2\lambda\gamma\|\widetilde{A}_{(1)}^k-A_{(1)}^k\|_F^2 + 2\gamma\langle (\widetilde{A}_{(1)}^k-A_{(1)}^k), (\widetilde{A}_{(1)}^k-A_{(1)}^k)F^k(F^k)^T \rangle -\gamma^2\|(\widetilde{A}_{(1)}^k-A_{(1)}^k)F^k\|_F^2 \\
      &=& 2\lambda\gamma\|\widetilde{\A}^k-\A^k\|_F^2 + (2\gamma -\gamma^2)\|(\widetilde{A}_{(1)}^k-A_{(1)}^k)F^k\|_F^2 \\
      &\ge& 2\lambda\gamma\|\widetilde{\A}^k-\A^k\|_F^2 \\
      &=& \frac{2\lambda}{\gamma}\|\A^{k+1}-\A^k\|_F^2,
  \end{eqnarray*}
  where the fifth equality holds because of \eqref{aaa}.

  Finally, in a similar way, we obtain
  \begin{eqnarray*}
    f(\X^{k+1},\A^{k+1},\B^k,\C^k)-f(\X^{k+1},\A^{k+1},\B^{k+1},\C^k)
      &\ge& \frac{2\lambda}{\gamma}\|\B^{k+1}-\B^k\|_F^2, \\
    f(\X^{k+1},\A^{k+1},\B^{k+1},\C^k)-f(\X^{k+1},\A^{k+1},\B^{k+1},\C^{k+1})
      &\ge& \frac{2\lambda}{\gamma}\|\C^{k+1}-\C^k\|_F^2.
  \end{eqnarray*}
  The lemma follows by summing the above four equalities.
\qed

\begin{lemma}\label{Lem-5.4}
  Let $\{\y^k\}_{k=0,1,2,\dots}$ be a sequence of iterates generated by Algorithm \ref{ALSb}. Then
  \begin{equation*}
    \sum_{k=0}^\infty \|\y^k-\y^{k+1}\|^2 < \infty
  \end{equation*}
  and
  \begin{equation*}
    \lim_{k\to\infty} \y^k-\y^{k+1}=0.
  \end{equation*}
\end{lemma}
{\bf Proof}
  From Lemma \ref{Lem-5.3}, we know
  \begin{eqnarray*}
    \left\|\y^k-\y^{k+1}\right\|^2 \le \frac{\gamma}{2\lambda}\left(f(\y^k)-f(\y^{k+1})\right),
  \end{eqnarray*}
  for $k=0,1,2,\dots$. By summarizing all $k$, we have
  \begin{eqnarray*}
    \sum_{k=0}^\infty \left\|\y^k-\y^{k+1}\right\|^2
      &\le& \frac{\gamma}{2\lambda}\sum_{k=0}^\infty \left(f(\y^k) - f(\y^{k+1})\right) \\
      &\le& \frac{\gamma}{2\lambda} f(\y^0) < \infty.
  \end{eqnarray*}
  where the second inequality holds because $f$ is always nonnegative. Hence, $\|\y^k-\y^{k+1}\|^2\to0$ and hence $\|\y^k-\y^{k+1}\|\to0$. The lemma is proved.
\qed

\begin{theorem}\label{Thm-5.5}
  Suppose that the infinite sequence of iterates $\{\y^k\}$ generated by Algorithm \ref{ALSb} is bounded. Then, every limit point of $\{\y^k\}$ is a stationary point.
\end{theorem}
{\bf Proof}
  Since $\{\y^k\}=\{(\X^k,\A^k,\B^k,\C^k)\}$ is bounded, it must have a subsequence $\{\y^{k_i}\}=\{(\X^{k_i},\A^{k_i},\B^{k_i},\C^{k_i})\}$ that converges to a limit point $\y^{\infty}=(\X^{\infty},\A^{\infty},\B^{\infty},\C^{\infty})$. Furthermore, the subsequence $\{\y^{k_i+1}\}=\{(\X^{k_i+1},\A^{k_i+1},\B^{k_i+1},\C^{k_i+1})\}$ also converges to the limit point $\y^{\infty}$ by Lemma \ref{Lem-5.4}.

  Because $\X^{k_i}-\X^{k_i+1}\to0$, we get $\widetilde{\X}^{k_i}-\X^{k_i}=0$ as $i\to\infty$. It yields from \eqref{xxx} that $\left(I-P^T(PP^T)^{-1}P\right)\mathrm{vec}(\A^{k_i}\B^{k_i}\C^{k_i}-\X^{k_i})\to0$ as $i\to\infty$. That is to say $\left(I-P^T(PP^T)^{-1}P\right)\mathrm{vec}(\A^\infty\B^\infty\C^\infty-\X^\infty)=0$.

  By $\A^{k_i}-\A^{k_i+1}\to0$, we have $\widetilde{\A}^{k_i}-\A^{k_i}\to0$. Because $\B^{k_i}\to\B^{\infty}$ and $\C^{k_i}\to\C^{\infty}$ as $i\to\infty$, the subsequence $\{F^{k_i}\}$ converges to $F^{\infty}$ that is bounded above. It yields from \eqref{aaa} that
  \begin{eqnarray*}
    (T^\infty_{(1)}-A_{(1)}^{\infty} F^{\infty}){F^{\infty}}^T
    &=& \lim_{i\to\infty} (T^{k_i+1}_{(1)}-A_{(1)}^{k_i}F^{k_i}){F^{k_i}}^T \\
    &=& \lim_{i\to\infty} (\widetilde{\A}^{k_i}-\A^{k_i})\left(F^{k_i}\left(F^{k_i}\right)^T+\lambda I_{r^2}\right) \\
    &=& 0.
  \end{eqnarray*}

  Finally, by a similar discussion, we know that $(X_{(2)}^{\infty}-B_{(2)}^{\infty}G^{\infty})(G^{\infty})^T = 0$ and $(X_{(3)}^{\infty}-C_{(3)}^{\infty}H^{\infty})(H^{\infty})^T = 0$.
  Hence, by Lemma \ref{Lem-5.1}, $\y^{\infty}=(\X^\infty,\A^\infty,\B^\infty,\C^\infty)$ is a stationary point of \eqref{CA-1}.
\qed

Theorem \ref{Thm-5.5} shows that every limit point of iterates generated by Algorithm \ref{ALSb} is a stationary point. Next, using the Kurdyka-{\L}ojasiewicz (KL) property \cite{Lojasiewicz-1963,Bolte-etal-2007,XY13}, we prove that the sequence of iterates from Algorithm \ref{ALSb} converges to a stationary point. The analysis in the remainder of this section is based on the outline of \cite{AttouchBolte-2009,Attouch-etal-2010}. Since $f(\y)+\delta_{\Omega}(\y)$ is a semi-algebraic function, where $\delta_{\Omega}(\cdot)$ is an indicator function defined on the affine manifold $\Omega$, the following KL inequality holds.

\begin{definition}[Kurdyka-{\L}ojasiewicz (KL) property]\label{Def-5.6}
  Let $U\in\Re^n$ be an open set and $f\,:\,U\to\Re$ be a semi-algebraic function. For every critical point $\y^*\in U$ of $f$, there is a neighborhood $V\subseteq U$ of $\y^*$, an exponent $\theta\in[\frac{1}{2},1)$ and a positive constant $\mu$ such that
  \begin{equation*}
    |f(\y)-f(\y^*)|^{\theta} \leq \mu\|\Pi_{\Omega}(\nabla f(\y))\|,
    \qquad \forall\,\y\in V,
  \end{equation*}
  where $\Pi_{\Omega}(\nabla f(\y))$ is defined by \eqref{pgy}.
\end{definition}

Next, we give a lower bound on the progress made by one iteration.

\begin{lemma}\label{Lem-5.7}
  Suppose that the infinite sequence $\{\y^k\}$ generated by Algorithm \ref{ALSb} is bounded.
  Then, there is a positive constant $\varsigma$ such that
  \begin{equation*}
    \|\y^{k}-\y^{k+1}\| \ge \varsigma\|\Pi_{\Omega}(\nabla f(\y^k))\|.
  \end{equation*}
\end{lemma}
{\bf Proof}
  From \eqref{xxx} and $\X^{k+1}-\X^k=\gamma(\widetilde{\X}^k-\X^k)$, it yields that
  \begin{eqnarray*}
    \|2(I-P^T(PP^T)^{-1}P^T)\vt{\A^k\B^k\C^k-\X^k}\|^2 &=& 4(1+\lambda)^2\|\vt{\widetilde{\X}^k-\X^k}\|^2 \\
      &=& \frac{4(1+\lambda)^2}{\gamma^2}\|\vt{\X^{k+1}-\X^k}\|^2.
  \end{eqnarray*}

  Since $\y^k=(\X^k,\A^k,\B^k,\C^k)$ is bounded, by \eqref{Fk-def}, \eqref{Gk-def}, and \eqref{Hk-def}, we could assume that
  \begin{equation*}
    \|F^k\|_F\le \kappa, \quad  \|G^k\|_F\le \kappa, \quad  \|H^k\|_F\le \kappa,
  \end{equation*}
  where $\kappa$ is a positive constant.
  By \eqref{aaa} and $A_{(1)}^{k+1}-A_{(1)}^k=\gamma(\widetilde{A}_{(1)}^k-A_{(1)}^k)$, we have
  \begin{eqnarray*}
    2(X_{(1)}^k-A_{(1)}^kF^k)(F^k)^T &=& 2(X_{(1)}^{k+1}-A_{(1)}^kF^k)(F^k)^T - 2(X_{(1)}^{k+1}-X_{(1)}^k)(F^k)^T \\
      &=& 2(\widetilde{A}_{(1)}^k-A_{(1)}^k)(F^k(F^k)^T+\lambda I_{r^2}) -2(X_{(1)}^{k+1}-X_{(1)}^k)(F^k)^T \\
      &=& \frac{2}{\gamma}(A_{(1)}^{k+1}-A_{(1)}^k)(F^k(F^k)^T+\lambda I_{r^2}) -2(X_{(1)}^{k+1}-X_{(1)}^k)(F^k)^T.
  \end{eqnarray*}
  Hence,
  \begin{eqnarray*}
    \lefteqn{ \|2(X_{(1)}^k-A_{(1)}^kF^k)(F^k)^T\|_F^2 }\\
      &\le& 2\left\|\frac{2}{\gamma}(A_{(1)}^{k+1}-A_{(1)}^k)(F^k(F^k)^T+\lambda I_{r^2})\right\|_F^2 + 2\left\|2(X_{(1)}^{k+1}-X_{(1)}^k)(F^k)^T\right\|_F^2 \\
      &\le& \frac{8(\kappa^2+\lambda r^2)^2}{\gamma^2}\|A_{(1)}^{k+1}-A_{(1)}^k\|_F^2 + 8\kappa^2\|X_{(1)}^{k+1}-X_{(1)}^k\|_F^2.
  \end{eqnarray*}

  Similarly, we can establish
  \begin{eqnarray*}
    \|2(X_{(2)}^k-B_{(2)}^kG^k)(G^k)^T\|_F^2 &=&
      \frac{8(\kappa^2+\lambda r^2)^2}{\gamma^2}\|B_{(2)}^{k+1}-B_{(2)}^k\|_F^2 + 8\kappa^2\|X_{(2)}^{k+1}-X_{(2)}^k\|_F^2, \\
    \|2(X_{(3)}^k-C_{(3)}^kH^k)(H^k)^T\|_F^2 &=&
      \frac{8(\kappa^2+\lambda r^2)^2}{\gamma^2}\|C_{(3)}^{k+1}-C_{(3)}^k\|_F^2 + 8\kappa^2\|X_{(3)}^{k+1}-X_{(3)}^k\|_F^2,
  \end{eqnarray*}

  In sum, we have
  \begin{eqnarray*}
    \lefteqn{ \|\Pi_{\Omega}(\nabla f(\y^k))\|^2 = \|2(I-P^T(PP^T)^{-1}P^T)\vt{\A^k\B^k\C^k-\X^k}\|^2 } \\
      &&~{} + \|2(X_{(1)}^k-A_{(1)}^kF^k)(F^k)^T\|_F^2 + \|2(X_{(2)}^k-B_{(2)}^kG^k)(G^k)^T\|_F^2 + \|2(X_{(3)}^k-C_{(3)}^kH^k)(H^k)^T\|_F^2 \\
      &\le& \left(\frac{4(1+\lambda)^2}{\gamma^2}+24\kappa^2\right)\|\X^{k+1}-\X^k\|_F^2 \\
      &&~{}+ \frac{8(\kappa^2+\lambda r^2)^2}{\gamma^2}\left( \|\A^{k+1}-\A^k\|_F^2 + \|\B^{k+1}-\B^k\|_F^2 + \|\C^{k+1}-\C^k\|_F^2\right) \\
      &\le& \max\left\{\frac{4(1+\lambda)^2}{\gamma^2}+24\kappa^2, \frac{8(\kappa^2+\lambda r^2)^2}{\gamma^2}\right\}\|\y^{k+1}-\y^k\|^2.
  \end{eqnarray*}
  This lemma is valid when we set $\varsigma^{-2}:=\max\left\{\frac{4(1+\lambda)^2}{\gamma^2}+24\kappa^2, \frac{8(\kappa^2+\lambda r^2)^2}{\gamma^2}\right\}$.
\qed

\begin{lemma}\label{Lem-5.8}
  Let $\y^*$ be one of the limiting points of $\{\y^k\}$. Assume that $\y^0$ satisfies
  $\y^0\in B(\y^*,\rho) \subseteq V$ where
  \begin{equation}\label{eq:4.2.3}
    \rho > \frac{\gamma\mu}{2\lambda\varsigma(1-\theta)}|f(\y^0)-f(\y^*)|^{1-\theta}
           + \|\y^0-\y^*\|.
  \end{equation}
  Then, we have the following assertions:
  \begin{equation}\label{neighbor}
    \y^k \in \mathbb{B}(\y^*,\rho),
    \qquad \forall\,k=0,1,\ldots
  \end{equation}
  and
  \begin{equation}\label{ixixixx}
    \sum_{k=0}^{\infty}\|\y^k-\y^{k+1}\|
    \leq \frac{\gamma\mu}{2\lambda\varsigma(1-\theta)}|f(\y^0)-f(\y^*)|^{1-\theta}.
  \end{equation}
\end{lemma}
{\bf Proof}
  We prove (\ref{neighbor}) by induction. Obviously, $\y^0 \in B(\y^*,\rho)$ when $k=0$.
  Second, we assume there is an integer $K$ such that
  \begin{equation*}
    \y^k \in B(\y^*,\rho), \qquad \forall\, 0\leq k\leq K,
  \end{equation*}
  which means that the KL property holds at these points. Now, we are going to prove that $\y^{K+1} \in B(\y^*,\rho)$.

  We define a scalar function
  \begin{equation}\label{eq:4.2.4}
    \varphi(\alpha) := \frac{1}{1-\theta}|\alpha-f(\y^*)|^{1-\theta}.
  \end{equation}
  It is easy to see that $\varphi(\cdot)$ is a concave function and $\varphi'(\alpha)=|\alpha-f(\y^*)|^{-\theta}$ if $\alpha \geq f(\y^*)$. Then, for $0\leq k \leq K$, it yields that
  \begin{eqnarray*}
    \varphi(f(\y^k))-\varphi(f(\y^{k+1}))
      &\geq& \varphi'(f(\y^k))(f(\y^k)-f(\y^{k+1})) \\
      &\geq& \frac{1}{|f(\y^k)-f(\y^*)|^{\theta}}\frac{2\lambda}{\gamma}\|\y^k-\y^{k+1}\|^2 \qquad[\text{Lemma \ref{Lem-5.3}}]\\
      &\geq& \frac{2\lambda}{\gamma\mu}\frac{\|\y^k-\y^{k+1}\|^2}{\|\Pi_{\Omega}(\nabla f(\y^k))\|}
                                              \qquad[\text{KL property}]\\
      &\geq& \frac{2\lambda\varsigma}{\gamma\mu}\frac{\|\y^k-\y^{k+1}\|^2}{\|\y^k-\y^{k+1}\|}
                                              \qquad[\text{Lemma \ref{Lem-5.7}}]\\
      &=& \frac{2\lambda\varsigma}{\gamma\mu}\|\y^k-\y^{k+1}\|.
  \end{eqnarray*}
  By summating $k$ from $0$ to $K$, we have
  \begin{eqnarray}
    \sum_{k=0}^K\|\y^k-\y^{k+1}\|
      &\leq& \frac{\gamma\mu}{2\lambda\varsigma} \sum_{k=0}^K [\varphi(f(\y^k))-\varphi(f(\y^{k+1}))] \nonumber\\
      &=& \frac{\gamma\mu}{2\lambda\varsigma} [\varphi(f(\y^0))-\varphi(f(\y^{K+1}))] \nonumber\\
      &\leq& \frac{\gamma\mu}{2\lambda\varsigma} \varphi(f(\y^0)). \label{eq:4.6}
  \end{eqnarray}
  Hence, it follows from (\ref{eq:4.6}) and (\ref{eq:4.2.3}) that
  \begin{equation*}
    \|\y^{K+1}-\y^*\|
      \leq \sum_{k=0}^K\|\y^{k+1}-\y^k\| + \|\y^0-\y^*\|
      \leq \frac{\gamma\mu}{2\lambda\varsigma} \varphi(f(\y^0)) + \|\y^0-\y^*\|
      < \rho
  \end{equation*}
  which means (\ref{neighbor}) holds. Moreover, we obtain (\ref{ixixixx}) by letting $K\to\infty$ in (\ref{eq:4.6}) and using (\ref{eq:4.2.4}).
\qed

\begin{theorem}\label{Thm-5.9}
  Assume that Algorithm \ref{ALSb} produces a bounded sequence $\{\y^k\}$. Then,
  \begin{equation*}
    \sum_{k=0}^{\infty} \|\y^{k+1}-\y^k\| < +\infty,
  \end{equation*}
  which implies that the entire sequence $\{\y^k\}$ converges.
\end{theorem}
{\bf Proof}
  Because $\{\y^k\}$ is bounded, it must have a limit point $\y^*$ and there is an index $k_0$ such that $\y^{k_0}\in B(\y^*,\rho)$.
  If we regard $\y^{k_0}$ as an initial point, Lemma \ref{Lem-5.8} holds. The entire sequence $\{\y^k\}$ satisfies \eqref{ixixixx}. Hence, this lemma is proved.
\qed

Finally, by consulting \cite{AttouchBolte-2009}, we give a result on the local convergence rate.

\begin{theorem}\label{Thm-5.10}
  Assume that Algorithm \ref{ALSb} produces a bounded sequence $\{\y^k\}$.  \\
  (1) If $\theta = \frac{1}{2}$, there exist $\eta>0$ and $\nu\in[0,1)$ such that
      \begin{equation*}
        \|\y^k - \y^*\| \leq \eta\nu^k,
      \end{equation*}
      which means that the sequences of iterates converges R-linearly. \\
  (2) If $\theta\in(\frac{1}{2},1)$, there exist $\eta>0$ such that
      \begin{equation*}
        \|\y^k - \y^*\| \leq \eta k^{-\frac{1-\theta}{2\theta-1}}.
      \end{equation*}
\end{theorem}
{\bf Proof}
  Without loss of generality, we assume that $\y^0\in B(\y^*,\rho)$. Let
  \begin{equation}\label{aaaaaa}
    \Delta_k := \sum_{i=k}^{\infty}\|\y^{i}-\y^{i+1}\| \geq \|\y^{k}-\y^*\|.
  \end{equation}
  From Lemma \ref{Lem-5.8}, we have
  \begin{eqnarray}
    \Delta_k &\leq& \frac{\gamma\mu}{2\lambda\varsigma(1-\theta)}|f(\y^k)-f(\y^*)|^{1-\theta} \nonumber\\
      &=& \frac{\gamma\mu}{2\lambda\varsigma(1-\theta)}\left[|f(\y^k)-f(\y^*)|^\theta\right]^{\frac{1-
      \theta}{\theta}} \nonumber\\
      &\leq& \frac{\gamma\mu}{2\lambda\varsigma(1-\theta)}\mu^{\frac{1-\theta}{\theta}}
             \|\Pi_{\Omega}(\nabla f(\y^k))\|^{\frac{1-\theta}{\theta}}
             \qquad [\text{KL property}]\nonumber\\
      &\leq& \frac{\gamma\mu}{2\lambda\varsigma(1-\theta)}\left(\frac{\mu}{\varsigma}\right)^{\frac{1-\theta}{\theta}}
             \|\y^k-\y^{k+1}\|^{\frac{1-\theta}{\theta}}
             \qquad [\text{Lemma \ref{Lem-5.7}}]\nonumber\\
      &=& \frac{\gamma\mu^{1/\theta}}{2\lambda\varsigma^{1/\theta}(1-\theta)}
             \|\y^k-\y^{k+1}\|^{\frac{1-\theta}{\theta}}. \label{11111}
  \end{eqnarray}

  (1) In the case $\theta = \frac{1}{2}$, we have $\frac{1-\theta}{\theta} = 1$ immediately. Then, the inequality (\ref{11111}) gives
  \begin{equation*}
    \Delta_k \leq \frac{\gamma\mu^{1/\theta}}{2\lambda\varsigma^{1/\theta}(1-\theta)} (\Delta_k-\Delta_{k+1}),
  \end{equation*}
  which means that
  \begin{equation}\label{bbbbbb}
    \Delta_{k+1} \leq \frac{{\gamma\mu^{1/\theta}}-{2\lambda\varsigma^{1/\theta}(1-\theta)}}{{\gamma\mu^{1/\theta}}} \Delta_k.
  \end{equation}
  Let $\nu:=\frac{{\gamma\mu^{1/\theta}}-{2\lambda\varsigma^{1/\theta}(1-\theta)}}{{\gamma\mu^{1/\theta}}}$.
  From (\ref{aaaaaa}) and (\ref{bbbbbb}), we know
  $\|\y^{k}-\y^*\| \leq \Delta_k \leq \nu\Delta_{k-1} \leq \cdots \leq \nu^k\Delta_0$,
  where $\Delta_0$ is finite by Theorem \ref{Thm-5.9}. Hence, assertion (1) is valid by taking $\eta:=\Delta_0$.

  (2) Let $\chi^{\frac{1-\theta}{\theta}} := \frac{\gamma\mu^{1/\theta}}{2\lambda\varsigma^{1/\theta}(1-\theta)}$.
  It yields from (\ref{11111}) that
  \begin{equation*}
    \Delta_k^{\frac{\theta}{1-\theta}} \leq \chi(\Delta_k-\Delta_{k+1}).
  \end{equation*}
  We define $h(\alpha):=\alpha^{-\frac{\theta}{1-\theta}}$. Obviously, $h(s)$ is monotonically decreasing. Then,
  \begin{eqnarray*}
    \frac{1}{\chi} &\leq& h(\Delta_k)(\Delta_k-\Delta_{k+1})\\
      &=& \int_{\Delta_{k+1}}^{\Delta_k}h(\Delta_k) \mathrm{d}\alpha \\
      &\leq& \int_{\Delta_{k+1}}^{\Delta_k}h(\alpha) \mathrm{d}\alpha \\
      &=& -\frac{1-\theta}{2\theta-1}
        \left(\Delta_k^{-\frac{2\theta-1}{1-\theta}}-\Delta_{k+1}^{-\frac{2\theta-1}{1-\theta}}\right).
  \end{eqnarray*}
  We denote $\vartheta:= -\frac{2\theta-1}{1-\theta}<0$ since $\theta\in(\frac{1}{2},1)$.
  Then, in follows from the above inequality that
  \begin{equation*}
    \Delta_{k+1}^\vartheta - \Delta_k^\vartheta \geq \frac{-\vartheta}{\chi}=:\varpi >0,
  \end{equation*}
  which gives
  \begin{equation*}
    \Delta_k \leq [\Delta_0^{\vartheta} + k\varpi]^{\frac{1}{\vartheta}}
      \leq (k\varpi)^{\frac{1}{\vartheta}}.
  \end{equation*}
  We obtain the assertion (2) by taking $\eta:=\varpi^{\frac{1}{\vartheta}}$.
\qed

\section{Numerical Tests}

In this section, we are going to compare the triple decomposition tensor recovery model \eqref{Ten-Recov-1} with the CP decomposition tensor recovery model and the Tucker decomposition tensor recovery model. Let $A\in\Re^{n_1\times r}$, $B\in\Re^{n_2\times r}$, and $C\in\Re^{n_3\times r}$. The CP tensor $[[A,B,C]]\in\Re^{n_1\times n_2\times n_3}$ has entries
\begin{equation*}
  [[A,B,C]]_{ijt} = \sum_{p=1}^r a_{ip}b_{jp}c_{tp}.
\end{equation*}
In addition, let $\D\in\Re^{r\times r\times r}$ be a core tensor. The Tucker tensor $[[\D; A,B,C]]$ has entries
\begin{equation*}
  [[\D; A,B,C]]_{ijt} = \sum_{p,q,s=1}^r  a_{ip}b_{jq}c_{ts}d_{pqs}.
\end{equation*}
Then, the CP tensor recovery model and the Tucker tensor recovery model could be represented by
\begin{equation*}
\begin{aligned}
  \min~ & \left\| [[A,B,C]] - \T \right\|_F^2 \\
  \st~  & \mathbb{P}(\T)=\dd
\end{aligned}
\end{equation*}
and
\begin{equation*}
\begin{aligned}
  \min~ & \| [[\D; A,B,C]] - \T \|_F^2 \\
  \st~  & \mathbb{P}(\T)=\dd,
\end{aligned}
\end{equation*}
respectively. These two models are solved by variants of MALS in Algorithm \ref{ALSb}.

\subsection{ORL Face Data}

\begin{figure}[tbh]
  \centering
  \includegraphics[width=\textwidth]{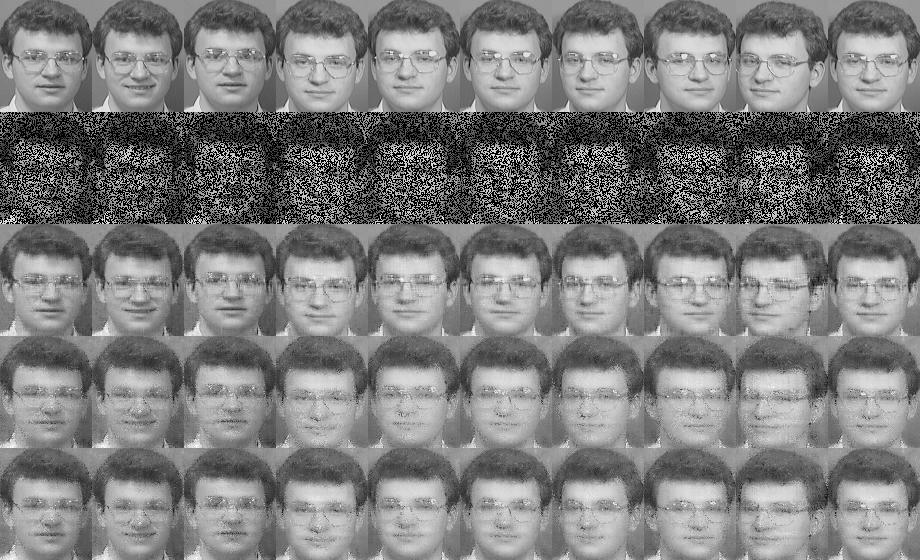}\\
  \caption{Original images are illustrated in the first row. Samples of 50 percent pixels are illustrated in the second row. The third to last rows report the recovered images by the proposed method, CP tensor recovery, and Tucker tensor recovery, respectively.}\label{ORL-C}
\end{figure}

Next, we apply the triple decomposition tensor recovery method, the CP decomposition tensor recovery method, and the Tucker decomposition tensor recovery method for the ORL dataset of faces. Original images of a person are illustrated in the first line of Figure \ref{ORL-C}. We sample fifty percent of pixels of these images as shown in the second line of Figure \ref{ORL-C}.

\begin{table}[htb]
  \centering
\begin{tabular}{c|ccccccc}
  \hline
  Methods & $r=1$ & $r=2$ & $r=3$ & $r=4$ & $r=5$ & $r=6$ & $r=7$ \\
  \hline
  New method    &  0.1496  &  0.1020  &  0.0850  &  0.0766  &  0.0712  &  0.0675  &  0.0667 \\
  CP recovery   &  0.1496  &  0.1127  &  0.1067  &  0.1021  &  0.0984  &  0.0940  &  0.0908 \\
  Tucker recov. &  0.1496  &  0.1127  &  0.1054  &  0.0993  &  0.0950  &  0.0894  &  0.0851 \\
  \hline
\end{tabular}
  \caption{Relative error of the recovered tensors.}\label{ORL-D}
\end{table}

Once a tensor $\T_{rec}$ is recovered, we calculate the relative error
\begin{equation*}
  \mathrm{RE} = \frac{\|\T_{rec}-\T_{truth}\|_F}{\|\T_{truth}\|_F}.
\end{equation*}
By varying rank $r$ from one to seven, the relative error of recovered tensor by the triple decomposition tensor recovery method, the CP decomposition tensor recovery method, and the Tucker decomposition tensor recovery method are reported in Table \ref{ORL-D}. Obviously, when the rank is one, all the models are equivalent. As the rank increases, the relative error of the recovered tensor by each method decreases. It is easy to see that the relative error corresponding to the proposed triple tensor recovery method decreases quickly. Hence, the new method performs better than the CP decomposition tensor recovery method and the Tucker decomposition tensor recovery method. We take rank $r=7$ for instance. The recovered tensors by the triple tensor recovery method, the CP decomposition tensor recovery method, and the Tucker
decomposition tensor recovery method are illustrated in lines 3--5 of Figure \ref{ORL-C}. Clearly, the quality of images recovered by proposed triple decomposition tensor recovery method is better than the others.

\subsection{McGill Calibrated Colour Image Data}

\begin{figure}[tbh]
  \centering
  \includegraphics[width=\textwidth]{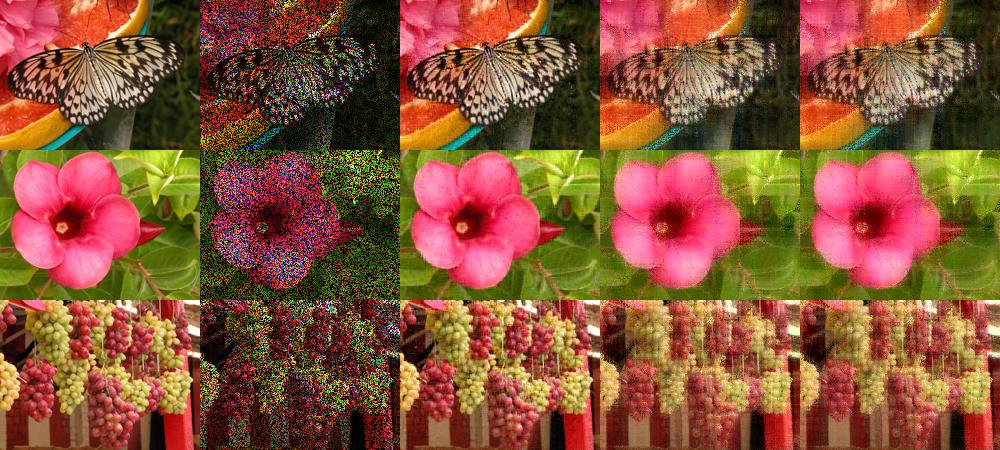}\\
  \caption{Original images are illustrated in the first column. Samples of 50 percent pixels are illustrated in the second column. The third to last columns report the recovered images by the proposed method, CP tensor recovery, and Tucker tensor recovery, respectively.}\label{McGill-B}
\end{figure}

We now investigate the McGill Calibrated Colour Image Data  \cite{CSZZC19, OK04}. We choose three colour images: butterfly, flower, and grape, which are illustrated in the first column of Figure \ref{McGill-B}. By resizing the colour image, we get a $150$-by-$200$-by-$3$ tensor. We randomly choose fifty percent entries of the colour image tensor. Tensors with missing entries are shown in the second column of Figure \ref{McGill-B}. We choose rank $r=7$, the proposed triple decomposition tensor recovery method generated an estimated tensor with relative error $0.1076$. Estimated color images are illustrated in the third column. Colour images estimated by the CP decomposition tensor recovery method and the Tucker decomposition tensor recovery method are shown in the fourth and the last columns with relative error $0.1921$ and $0.1901$, respectively. Obviously, the proposed triple tensor recovery method generates unambiguous colour images.

\section{Concluding Remarks}

In this paper, we introduce triple decomposition and triple rank for third order tensors.   The triple rank of a third order tensor is not greater than the CP rank and the middle value of the Tucker rank, is strictly less than the CP rank with a substantial probability, and is strictly less than the middle value of the Tucker rank for an essential class of examples. This indicates that practical data can be approximated by low rank triple decomposition as long as it can be approximated by low rank CP or Tucker decomposition.    We confirm this theoretical discovery numerically.  Numerical tests show that third order tensor data from practical applications such as internet traffic and video image are of low triple ranks. Finally, we have considered an application of triple decomposition to tensor recovery. Given its simplicity and practical applicability, we conclude that further study on triple decomposition is worth being conducted.

\bigskip

{\bf Acknowledgment}  We are thankful to Haibin Chen and Ziyan Luo for their comments, and to Qun Wang, who drew Figure 1 for us.

\end{document}